# LOCUS OF NON-REAL EIGENVALUES OF A CLASS OF LINEAR RELATIONS IN A KREIN SPACE

RYTIS JURŠĖNAS

ABSTRACT. It is a classical result that, if a maximal symmetric operator $T$ in a Krein space $\mathscr{H} = \mathscr{H}^-[\oplus]\mathscr{H}^+$ has the property $\mathscr{H}^- \subseteq \mathscr{D}_T$, then the imaginary part of its eigenvalue $\lambda$ from upper or lower half-plane is bounded by $|\mathrm{Im}\,\lambda| \leq 2\|TP^-\|$. We prove that in both half-planes $|\mathrm{Im}\,\lambda|$ never exceeds $t_0\|TP^-\|$ for some constant $t_0 \approx 1.84$. The result applies to a closed symmetric relation $T$ and carries on a suitable, most notably dissipative, extension.

## 1. INTRODUCTION

Let $T$ be a relation in a Krein space $(\mathscr{H}, [\cdot\,,\,\cdot])$, with the property
$(L)$ $\mathscr{D}_T$ contains a maximal uniformly negative subspace.

Stated otherwise, there is a canonical decomposition $\mathscr{H} = \mathscr{H}^-[\oplus]\mathscr{H}^+$ with negative subspace $\mathscr{H}^-$ contained in the domain $\mathscr{D}_T$.

With the opposite sign in an indefinite metric $[\cdot\,,\,\cdot]$ $T \in (L)$ agrees with [6, Definition 3.1.5], followed by [37] (and citation there); one finds examples and further discussion in [7, 9, 45]. Here we stick to the above-stated $(L)$, for it suits a (originally defined [39]) Pontryagin space with finite negative index.

If $T \in (L)$ is a maximal dissipative operator in $\mathscr{H}$, its resolvent set $\rho(T)$ contains those $\lambda \in \mathbb{C}^-$ with $\mathrm{Im}\,\lambda < -2\|TP^-\|$ (*cf.* [6, Theorem 2.2.9]). In particular this implies that $-2\|TP^-\| \leq \mathrm{Im}\,\lambda < 0$ for an eigenvalue $\lambda$ of $T$ from $\mathbb{C}^-$. If $T$ is maximal symmetric, an analogous estimate of $\mathrm{Im}\,\lambda$ from upper half-plane may exist too (*cf.* [6, Corollary 2.3.9]).

A maximal dissipative operator in $\mathscr{H}$ is necessarily densely defined (recall [6, Theorem 2.2.6], [8, Lemma 2.1]), the latter being a special case of a relation $T$ in $\mathscr{H}$ with the property
$(P)$ $\mathscr{D}_T + \mathscr{R}_T$ is dense in $\mathscr{H}$.

Alike a maximal dissipative extension of a closed dissipative operator is an operator, a closed canonical extension with nonempty resolvent set of a closed symmetric operator with properties $(L) \wedge (P)$ is an operator.

---

Given a closed symmetric relation $T \in (L)$ in $\mathscr{H}$ ($\mathscr{H}^- \subseteq \mathscr{D}_T$), the aim is to show that $|\mathrm{Im}\,\lambda|$ never reaches $2\|T_s P^-\|$; $T_s$ is the operator part and $\lambda$ an eigenvalue of $T$. Points of regular type are estimated too.

For a precise statement, the loci are defined as follows: $\Gamma_T$ is the set of those $\lambda \in \mathbb{C}_* = \mathbb{C} \setminus \mathbb{R}$ such that

$$|\mathrm{Im}\,\lambda| > \|T_s P^-\|\,, \qquad \left|\sigma_\lambda + \frac{\|P^- T P^-\|}{\lambda}\right| < \frac{2}{1 + \left(\frac{\|T_s P^-\|}{|\mathrm{Im}\,\lambda|}\right)^2}\,,$$

$$\sigma_\lambda = \begin{cases} 1, & \mathrm{Re}\,\lambda \geq 0\,, \\ -1, & \mathrm{Re}\,\lambda < 0\,. \end{cases}$$

The complement $C_T = \mathbb{C}_* \setminus \Gamma_T$.

In a classic example when a self-adjoint block operator $T$ is generated by a quadratic operator pencil ([6, Section 4.3.5], [34]) one has $T(\mathscr{H}^-) \subseteq \mathscr{H}^+$, and then $C_T = \{\lambda \mid 0 < |\mathrm{Im}\,\lambda| \leq \|TP^-\|\}$.

In general one estimates $C_T$ as follows.

**Proposition 1.1.** *If $t_0 \approx 1.83929$ denotes the (only one) real zero of the complex polynomial $t^3 - t^2 - t - 1$, then*

$$C_T \subseteq \{\lambda \mid 0 < |\mathrm{Im}\,\lambda| \leq t_0 \|T_s P^-\|\}\,.$$

*Proof.* This is because $(1+\rho^2)|\lambda \pm \|P^- T_s P^-\|| < 2|\lambda|$, with $\rho = m/|\mathrm{Im}\,\lambda|$ and $m = \|T_s P^-\|$, always holds if $(1+\rho^2)(|\lambda|+m) < 2|\lambda|$ or equivalently if $(1+\rho^2)m < (1-\rho^2)|\lambda|$; the latter in turn always holds if $(1+\rho^2)m < (1-\rho^2)|\mathrm{Im}\,\lambda|$ *i.e.* if $t = \rho^{-1} > 1$ solves $t^3 - t^2 - t - 1 > 0$. □

*Example* 1.2. Sample of $\Gamma_T$ (gray) for $\|T_s P^-\| = 1$, $\|P^- T_s P^-\| = 0.8$:

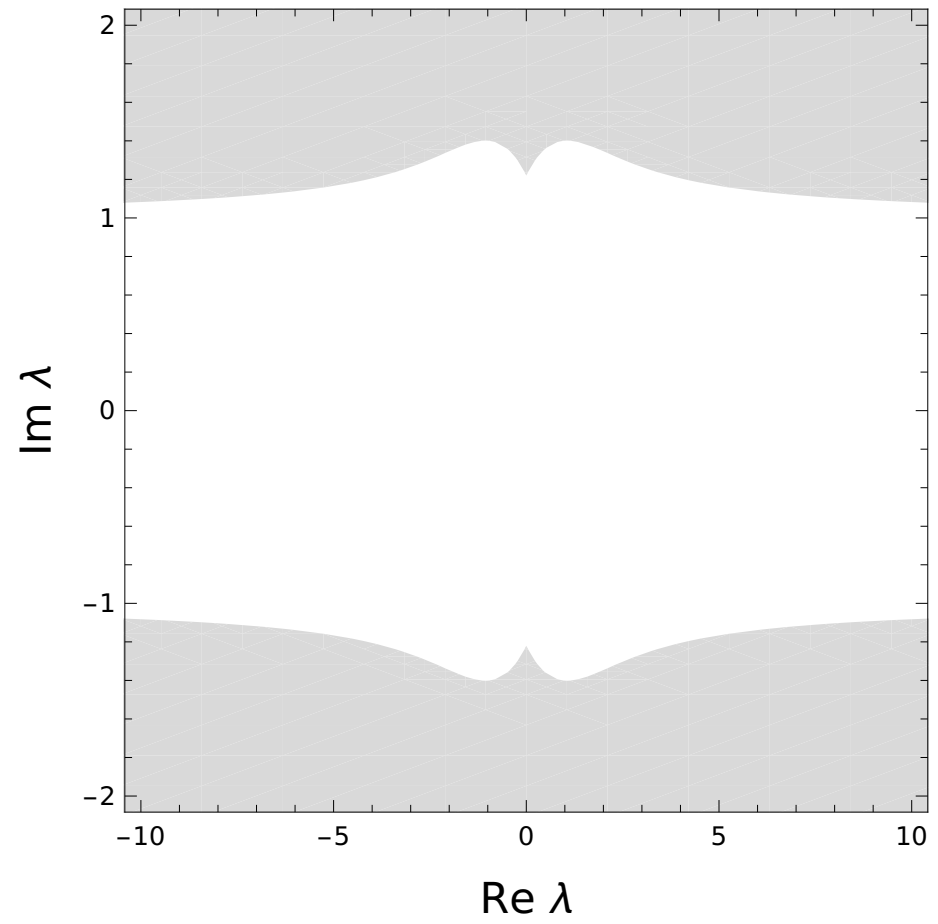

According to the next theorem, all non-real eigenvalues of a closed symmetric relation $T \in (L)$ are located within $C_T$.

**Theorem 1.3.** *Let $T \in (L)$ be a closed symmetric relation in a Krein space $\mathcal{H} = \mathcal{H}^-[\oplus]\mathcal{H}^+$, $\mathcal{H}^- \subseteq \mathcal{D}_T$, with canonical symmetry $J$, and let $T_0$ be a dissipative extension in $\mathcal{H}$. Then*

a) $\mathbb{C}^- \cap \sigma_p(T_0) \subseteq \mathbb{C}^- \cap C_T$.

b) *$T_0$ is an operator, if $T \in (P)$ and is an operator.*

c) *$\Delta_0^- \cap \Gamma_T \subseteq r(T_0)$, where $\Delta_0^-$ is the set of those $\lambda \in \mathbb{C}^-$ such that the lineals $\mathcal{H}^+ + \mathcal{R}_{JT_0+\lambda I}$ and $P^+(JT_0 - \lambda I)(JT_0 + \lambda I)^{-1}(\mathcal{H}^+)$ are closed. Particularly $\Delta_0^- = \mathbb{C}^-$ if $\dim \mathcal{H}^- < \infty$.*

As usual $\sigma_p(T_0)$ denotes the point spectrum and $r(T_0)$ is the field of regularity, and the set $\hat{\sigma}(T_0) = \mathbb{C} \setminus r(T_0)$ is known as the spectral core. Regardless of how one defines the continuous spectrum (*cf.* [6, Section 2.6], [20, Section VII.5.1] and [44, Section 2.2], [17, Section 1.7]), c) follows from the equality $\Delta_0^- \cap \hat{\sigma}(T_0) = \Delta_0^- \cap \sigma_p(T_0)$, provided one shows that $\lambda$ from $\Delta_0^-$ makes the lineal $\mathcal{R}_{T_0-\lambda I}$ into a subspace. As an example the latter equality is well-known to hold for a maximal $\pi$-dissipative operator $T_0$: in this case $\Delta_0^- = \mathbb{C}^-$ and moreover $\mathbb{C}^- \cap \hat{\sigma}(T_0)$ coincides with the part in $\mathbb{C}^-$ of the spectrum $\sigma(T_0)$, and hence $\mathbb{C}^- \cap \Gamma_T \subseteq \rho(T_0)$.

In Theorem 1.3 $T_0$ is necessarily a canonical extension, $T \subseteq T_0 \subseteq T^c$, by the Cauchy–Schwarz inequality for a non-negative lineal. A standard example of a maximal $J$-dissipative extension of a $J$-symmetric relation (or operator) is constructed from the Štraus extension [1, 2, 3, 4], [12, Proposition 1.6.8] (*cf.* [41, Theorem 4.7]) by means of multiplication by $J$.

As a special case $\mathbb{C}_* \cap \sigma_p(T_0) \subseteq C_T$ for a symmetric extension $T_0$; if $\mathcal{H}$ is a Pontryagin space, then also $\Gamma_T \subseteq r(T_0)$.

The theorem really becomes of interest only if a $J$-symmetric $T \neq T^c$, because the non-real spectrum of a $J$-self-adjoint relation $T \in (L)$ consists exclusively of normal eigenvalues and lies within the part $\{\lambda \mid 0 < |\operatorname{Im}\lambda| \leq \|T_s P^-\|\}$ of $C_T$, namely in $\{\lambda \mid 0 < |\operatorname{Im}\lambda| \leq \|P^+ T_s P^-\|\}$.

If in the theorem $\mathcal{H}^-$ (and therefore $\mathcal{H}^+$) is invariant for $T$, then $T$ is an operator without non-real eigenvalues.

It will be clear from the proofs, on the other hand, that "dissipative" in the theorem, however common in applications, can be replaced by a larger class—the only thing needed from dissipativeness is that $\mathbb{C}^- \subseteq r(JT_0)$, whatever the canonical symmetry $J$ is. The theorem, which we do not state explicitly, is modified appropriately.

Theorem 1.3 a) is shown in the next Section 2 by using the notion of an angular operator for a non-negative lineal in $\mathcal{H}$.

Theorem 1.3 b) is a particular result of discussion in Sections 3, 4.

A closed $J$-symmetric relation $T \in (L)$ with the property $(P)$ or

$(P')$ $\mathcal{D}_T + \operatorname{Ind} T$ is dense in $\mathcal{H}$

is studied in Sections 4, 5. The example of $T \in (L) \wedge (P')$ is a $J$-symmetric relation being 1-1 with the forbidden relation of a boundary relation; compared to [25] our construction deals with an extra $(L)$ (which is an attribute of a $J$-space) and removes the restriction imposed on defect numbers.

Theorem 1.3 c) follows from the two lemmas in Section 6. For a closed relation $T$ in a $J$-space, the idea is based on finding those $\lambda \in r(JT)$ for which $\mathscr{R}_{T-\lambda I}$ is closed. Notably, a nonzero $\lambda \in r(JT)$ is a $\Phi_+$-point ([21]) of a closed relation $T$ in a Pontryagin space.

Two rather classic, but viewed from the present perspective, spectral results for a $J$-self-adjoint relation are mentioned in the last Section 7.

## 2. LOCATION OF NON-REAL EIGENVALUES

2.1 We use standard terminology from Krein space [5, 6, 13, 29, 39] as well as Hilbert space [12, 17, 24, 28, 40] theory. Notation is essentially as in [6], with the (componentwise) sums of relations as in [12, 24].

A lineal is a linear set in a linear space. A relation $T$ in a linear space $\mathscr{H}$ is a lineal in $\mathscr{H}_\Gamma = \mathscr{H}^2 (= \mathscr{H} \times \mathscr{H})$ with domain $\mathscr{D}_T$, range $\mathscr{R}_T$, kernel $\operatorname{Ker} T$, also $\operatorname{Ker}_\lambda T = \operatorname{Ker}(T - \lambda I)$ ($\lambda \in \mathbb{C}$), and multivalued part (indefiniteness) $\operatorname{Ind} T$. Given a lineal $\mathscr{L} \subseteq \mathscr{H}$, the domain restriction $T \mid \mathscr{L} = T \mid \mathscr{L} \cap \mathscr{D}_T$, with the same meaning for its range $T(\mathscr{L})$.

An overbar refers to the closure if a Hilbert metric is assumed; whenever it matters an inner product is conjugate-linear in its first factor. An operator in a Hilbert space is a relation $T$ with $\operatorname{Ind} T$ trivial. A subspace of a Hilbert space is a closed lineal.

Let $(\mathscr{H}, [\cdot\,, \cdot])$ be a Krein space and $T$ a relation in $\mathscr{H}$.

**Lemma 2.1.** (*cf.* [31, Lemma 3.1]) *If $T$ is a closed relation in a Krein space and $\mathscr{L}$ is a subspace, then*

a) $(T(\mathscr{L}))^{[\perp]} = (T^c)^{-1}(\mathscr{L}^{[\perp]})$ *if $\mathscr{L} + \mathscr{D}_T$ is a subspace.*

b) $T(\mathscr{L})$ *is a subspace if such are $\mathscr{L} + \mathscr{D}_T$ and $\mathscr{L}^{[\perp]} + \mathscr{R}_{T^c}$.*

c) $\overline{T(\mathscr{L})} = T(\overline{\mathscr{L} + \operatorname{Ker} T})$ *if $\mathscr{L} + \mathscr{D}_T$ and $\mathscr{R}_T$ are subspaces.*

The orthogonal complement in $\mathscr{H}_\Gamma$ of a lineal $T (\subseteq \mathscr{H}_\Gamma)$ with respect to an indefinite metric $[(x_1, y_1)\,, (x_2, y_2)]_\Gamma = -\mathrm{i}([x_1\,, y_2] - [y_1\,, x_2])$ (*cf.* [6, Section 2.1]) is referred to as the adjoint (of $T$), [6, Definition 2.6.7], and denoted $T^c$. $T$ is dissipative (symmetric, self-adjoint) if it is a non-negative lineal (resp. a neutral lineal, a hyper-maximal neutral subspace) in the Krein space $(\mathscr{H}_\Gamma, [\cdot\,, \cdot]_\Gamma)$.

Given a canonical decomposition $\mathscr{H} = \mathscr{H}^- [\oplus] \mathscr{H}^+$ with canonical projections $P^\pm$ and canonical symmetry $J = P^+ - P^-$, a dissipative (symmetric, self-adjoint) relation $T$ is $J$-dissipative (resp. $J$-symmetric, $J$-self-adjoint) if explicit reference to $\mathscr{H}$ is omitted; $\mathscr{H}$ itself is a $J$-space. The same prefix "$J$-" is occasionally used for the adjoint $T^c$, as well as if the orthogonality of a lineal $\mathcal{L}(\subseteq \mathcal{H})$ is assumed with respect to $[\cdot\,,\,\cdot]$ (*i.e.* one makes a distinction between the $J$-orthogonal and the orthogonal complement, $\mathcal{L}^{[\perp]}$ and $\mathcal{L}^{\perp} = J(\mathcal{L}^{[\perp]})$ respectively). In a $\Pi_\kappa$-space (a Pontryagin space with finite number $\kappa$ of negative squares) "$J$-" is replaced by "$\pi$-".

Given a $J$-space, $\mathcal{H}_\Gamma$ is therefore a $J_\Gamma$-space, where the canonical symmetry $J_\Gamma(x, y) = (-\mathrm{i}Jy, \mathrm{i}Jx)$.

2.2 $\sigma_p(T) = \{\lambda \mid \operatorname{Ker}_\lambda T \neq \{0\}\}$ is the point spectrum of a relation $T$ in a linear space. If $T$ is closed in a Krein space $\mathcal{H}$, $r(T) = \{\lambda \mid \mathcal{R}_{T^c - \bar{\lambda}I} = \mathcal{H}\}$ is its field of regularity and $\rho(T) = \{\lambda \mid \mathcal{R}_{T^c - \bar{\lambda}I} = \mathcal{H} = \mathcal{R}_{T-\lambda I}\}$ is the resolvent set; *cf.* [36], [6, Section 2.6]. The spectrum $\sigma(T) = \mathbb{C} \setminus \rho(T)$.

The set $\mathbb{C}_* = \mathbb{C} \setminus \mathbb{R} = \mathbb{C}^+ \sqcup \mathbb{C}^-$ (disjoint $\cup$), $\mathbb{C}^\pm = \{\lambda \mid \pm \operatorname{Im}\lambda > 0\}$.

2.3 In a $J$-space $(\mathcal{H}, [\cdot\,,\,\cdot])$ the Hilbert metric $\langle x\,,\,y\rangle = [x\,,\,Jy]$ and the induced norm $\|x\| = \sqrt{[x\,,\,Jx]}$, and the operator (sup-)norm is denoted by the same $\|\cdot\|$.

Let $T$ be a closed relation in $\mathcal{H} = \mathcal{H}^-[\oplus]\mathcal{H}^+$; $T_s$ is its operator part. If $\mathcal{H}^- \subseteq \mathcal{D}_T$ (*i.e.* $T \in (L)$), $T_s \mid \mathcal{H}^-$ is bounded: $\|T_s \mid \mathcal{H}^-\| = \|T_s P^-\|$.

**Proposition 2.2.** (*cf.* [6, Corollary 2.3.9]) *Assume $T \in (L)$ is maximal symmetric in a Krein space with canonical decomposition $\mathcal{H} = \mathcal{H}^-[\oplus]\mathcal{H}^+$; $\mathcal{H}^- \subseteq \mathcal{D}_T$. Then either*
$\{\lambda \mid \operatorname{Im}\lambda > 2\|T_s P^-\|\} \subseteq \rho(T)$ *or*
$\{\lambda \mid \operatorname{Im}\lambda < -2\|T_s P^-\|\} \subseteq \rho(T)$.

*Proof.* This is evident from (*cf.* [6, Theorem 2.2.9])

$$
\begin{aligned}
T - \lambda I &\supseteq (TJ - \lambda I)(I + 2(TJ - \lambda I)^{-1} T P^-) \\
&\supseteq (TJ - \lambda I)(I + 2(TJ - \lambda I)^{-1} T_s P^-)
\end{aligned}
$$

for $\lambda \in \mathbb{C}$, $J = P^+ - P^-$: Let $(x, y) \in (TJ - \lambda I)(I + 2(TJ - \lambda I)^{-1} T_s P^-)$, *i.e.* $(\exists u)\,(\exists v)\,(Jx + Jv, y + \lambda x + \lambda v) \in T$, $(x - Jx, u) \in T_s$, $(Jv, u + \lambda v) \in T$; from here the domain and multivalued part of $(TJ - \lambda I)(I + 2(TJ - \lambda I)^{-1} T_s P^-)$ are given by $\{x \in \mathcal{D}_T \mid T_s P^- x \in \mathcal{R}_{TJ - \lambda I}\}$ and $\operatorname{Ind} T$ respectively, and moreover $(x, y) \in T - \lambda I$. □

*Remark* 2.3. In the proof and in some other proofs we use tacitly the next well-known criterion (*e.g.* [40, Proposition 2.1]): If $S \subseteq T$ are relations in a linear space, $S = T$ iff $\mathcal{D}_S = \mathcal{D}_T$ and $\operatorname{Ind} S = \operatorname{Ind} T$.

The proof allows one to formulate an analogous proposition for a maximal $J$-dissipative relation: $\{\lambda \mid \operatorname{Im}\lambda < -2\|T_s P^-\|\} \subseteq \rho(T)$.

2.4 In Proposition 2.2 $\mathbb{C}^\pm \cap \sigma_p(T) \subseteq \{\lambda \in \mathbb{C}^\pm \mid |\operatorname{Im}\lambda| \leq 2\|T_s P^-\|\}$ for a maximal symmetric $T \in (L)$. Before showing that actually $\mathbb{C}_* \cap \sigma_p(T) \subseteq C_T$ (that is Theorem 1.3 a)) we shall deal with some auxiliary propositions. Part of them we use later for other purposes.

**Lemma 2.4.** *Let $T \in (L)$ be a symmetric relation in a Krein space $(\mathcal{H},[\cdot,\cdot])$ with canonical decomposition $\mathcal{H}=\mathcal{H}^-[\oplus]\mathcal{H}^+$; $\mathcal{H}^- \subseteq \mathcal{D}_T$. Then*

a) *$P^-T^c\,|\,\mathcal{H}^- = P^-T\,|\,\mathcal{H}^-$ and is a bounded self-adjoint operator in a Hilbert space $(\mathcal{H}^-,-[\cdot,\cdot])$.*

b) *The adjoint of a closed relation $P^+T^c\,|\,\mathcal{H}^-\colon(\mathcal{H}^-,-[\cdot,\cdot])\to\mathcal{H}^+$ is the closure of an operator $-P^-\bar{T}\,|\,\mathcal{H}^+$ (which moreover coincides with the closure of an operator $-P^-T^c\,|\,\mathcal{H}^+$ if $\bar{T}\in(P')$).*

c) *The adjoint of a closed relation $P^+T^c\,|\,\mathcal{H}^+$ in $\mathcal{H}^+$ is a closed relation $P^+\bar{T}\,|\,\mathcal{H}^+$. (If $T$ is not necessarily symmetric in $\mathcal{H}$ and if one assumes $\mathcal{H}^-\subseteq\mathcal{D}_{T^c}$ instead of $\mathcal{H}^-\subseteq\mathcal{D}_T$, then the adjoint of a closed relation $P^+\bar{T}\,|\,\mathcal{H}^+$ in a Hilbert space $\mathcal{H}^+$ is $\overline{P^+T^c\,|\,\mathcal{H}^+}$.)*

d) *Let $T$ be closed. Regarding $P^-T^c\,|\,\mathcal{H}^+$ on $\mathcal{H}^+\cap\mathcal{D}_{T^c}$ as a densely defined Hilbert space operator from $\overline{\mathcal{H}^+\cap\mathcal{D}_{T^c}}(=\mathcal{H}^+\cap\bar{\mathcal{D}}_{T^c})$ to $(\mathcal{H}^-,-[\cdot,\cdot])$, the adjoint operator is given by $-P^+T_s\,|\,\mathcal{H}^-$.*

*Proof.* a) The equality is clear from

$$T^c\,|\,\mathcal{D}_T = T \,\hat{+}\, (\{0\}\times \operatorname{Ind} T^c)\,. \tag{2.1}$$

Using

$$P^-T\,|\,\mathcal{H}^- = (\mathcal{H}^-)^2 \cap (T \,\hat{+}\, (\{0\}\times\mathcal{H}^+))$$

the adjoint in $(\mathcal{H}^-,[\cdot,\cdot])$ of $P^-T\,|\,\mathcal{H}^-$, which coincides with the adjoint in $(\mathcal{H}^-,-[\cdot,\cdot])$, is a subspace of $\mathcal{H}_\Gamma$ that contains a subspace

$$(\mathcal{H}^-)^2\cap\overline{((\mathcal{H}^+)^2 \,\hat{+}\, (T^c\,|\,\mathcal{H}^-))} = (\mathcal{H}^-)^2\cap((\mathcal{H}^+)^2 \,\hat{\oplus}\, \overline{P^-T\,|\,\mathcal{H}^-})$$

*i.e.* it contains $P^-T\,|\,\mathcal{H}^-$. But $P^-T\,|\,\mathcal{H}^-$ is an everywhere defined operator in $\mathcal{H}^-$, so the above inclusions become equalities.

b) The relation

$$P^+T^c\,|\,\mathcal{H}^- = (\mathcal{H}^-\times\mathcal{H}^+)\cap(T^c \,\hat{+}\, (\{0\}\times\mathcal{H}^-))\,.$$

Because $T^c \,\hat{+}\, (\{0\}\times\mathcal{H}^-)$ is a closed relation in $\mathcal{H}$ (*i.e.* $\mathcal{D}_{\bar{T}}+\mathcal{H}^+=\mathcal{H}$), so is $P^+T^c\,|\,\mathcal{H}^-$. Regarding the latter as a Hilbert space relation $\mathcal{H}^-\to\mathcal{H}^+$ the corresponding adjoint, $(P^+T^c\,|\,\mathcal{H}^-)^*$, equals

$$-(\mathcal{H}^+\times\mathcal{H}^-)\cap\overline{((\mathcal{H}^-\times\mathcal{H}^+) \,\hat{+}\, (\bar{T}\,|\,\mathcal{H}^+))} = -\overline{P^-\bar{T}\,|\,\mathcal{H}^+}\,.$$

If moreover $\bar{T}\in(P')$ then $\operatorname{Ind}\bar{T}=\operatorname{Ind}T^c$ and then (see (2.1))

$$P^+T^c\,|\,\mathcal{H}^- = P^+\bar{T}\,|\,\mathcal{H}^- = (\mathcal{H}^-\times\mathcal{H}^+)\cap(\bar{T} \,\hat{+}\, (\{0\}\times\mathcal{H}^-))$$

so that then by arguing as above $(P^+T^c\,|\,\mathcal{H}^-)^* = -\overline{P^-T^c\,|\,\mathcal{H}^+}$.

c) Considering a closed relation

$$P^+T^c\,|\,\mathcal{H}^+ = (\mathcal{H}^+)^2\cap(T^c \,\hat{+}\, (\{0\}\times\mathcal{H}^-))$$

in $\mathcal{H}$ as a relation in $\mathcal{H}^+$, the adjoint is given by

$$(\mathcal{H}^+)^2\cap\overline{((\mathcal{H}^-)^2 \,\hat{+}\, (\bar{T}\,|\,\mathcal{H}^+))} = \overline{P^+\bar{T}\,|\,\mathcal{H}^+}\,.$$

Moreover $P^+\bar{T}\,|\,\mathcal{H}^+$ is a subspace in $\mathcal{H}_\Gamma$, since $\mathcal{D}_{T^c} + \mathcal{H}^+ = \mathcal{H}$.

If $\mathcal{H}^- \subseteq \mathcal{D}_{T^c}$ instead of $\mathcal{H}^- \subseteq \mathcal{D}_T$, then the adjoint of a closed relation $P^+\bar{T}\,|\,\mathcal{H}^+$ in $\mathcal{H}^+$ is $\overline{P^+T^c\,|\,\mathcal{H}^+}$: By arguing as above $\bar{T} \hat{+} (\{0\} \times \mathcal{H}^-)$, and hence $P^+\bar{T}\,|\,\mathcal{H}^+$, is a closed relation in $\mathcal{H}$. Viewing the latter in $\mathcal{H}^+$ the adjoint reads

$$(\mathcal{H}^+)^2 \cap \overline{((\mathcal{H}^-)^2 \hat{+} (T^c\,|\,\mathcal{H}^+))} = \overline{P^+T^c\,|\,\mathcal{H}^+}$$

which proves the assertion.

d) It suffices to show that the adjoint $(P^-T^c\,|\,\mathcal{H}^+)^* \supseteq -P^+T_s\,|\,\mathcal{H}^-$. First, because

$$P^-T^c\,|\,\mathcal{H}^+ = (\mathcal{H}^+ \times \mathcal{H}^-) \cap (T^c \hat{+} (\{0\} \times \mathcal{H}^+))$$

so that (with $\mathcal{H}_c^+ = \mathcal{H}^+ \cap \bar{\mathcal{D}}_{T^c}$)

$$-(P^-T^c\,|\,\mathcal{H}^+)^* \supseteq (\mathcal{H}^- \times \mathcal{H}_c^+) \cap P^+T\,|\,\mathcal{H}^-$$

and second, because

$$P^+T\,|\,\mathcal{H}^- = (P^+T_s\,|\,\mathcal{H}^-) \hat{\oplus} (\{0\} \times \operatorname{Ind} T)$$

so that $-(P^-T^c\,|\,\mathcal{H}^+)^* \supseteq P^+T_s\,|\,\mathcal{H}^-$. □

The next lemma is of interest only for $T \neq T^c$.

**Lemma 2.5.** *Let $T \in (L)$ be a closed symmetric relation in a Krein space $(\mathcal{H}, [\cdot\,,\cdot])$ with canonical decomposition $\mathcal{H} = \mathcal{H}^-[\oplus]\mathcal{H}^+$; $\mathcal{H}^- \subseteq \mathcal{D}_T$. To each fixed $\lambda \in \Gamma_T$ there corresponds a canonical decomposition $\mathcal{H} = \mathcal{H}_1^-[\oplus]\mathcal{H}_1^+$ with canonical symmetry $J_1$ such that* $\operatorname{Ker}_\lambda J_1 T^c \subseteq \mathcal{H}_1^+$.

*Proof.* Let $J$ be the canonical symmetry associated with the canonical decomposition $\mathcal{H} = \mathcal{H}^-[\oplus]\mathcal{H}^+$, $\mathcal{H}^- \subseteq \mathcal{D}_T$.

Then $\mathcal{N}_\lambda = \operatorname{Ker}_\lambda JT^c$ is uniformly positive if $|\operatorname{Im}\lambda| > m = \|T_sP^-\|$:

By Lemma 2.4 a) $\mathcal{H}^- \cap \mathcal{N}_\lambda = \{0\}$ (*i.e.* $\mathcal{N}_\lambda$ is of class $\mathcal{A}^+$ in [6, Section 1.8]), so $\mathcal{N}_\lambda = \{x^+ + K_\lambda x^+ \,|\, x^+ \in P^+(\mathcal{N}_\lambda)\}$; $K_\lambda = P^-(P^+\,|\,\mathcal{N}_\lambda)^{-1}$ is the restriction to $\mathcal{D}_{K_\lambda} = P^+(\mathcal{N}_\lambda)$ of the operator

$$\tilde{K}_\lambda = -(T^{--} + \lambda I^-)^{-1}(P^-T^c\,|\,\mathcal{H}^+)\,, \quad \mathcal{D}_{\tilde{K}_\lambda} = \mathcal{H}^+ \cap \mathcal{D}_{T^c}$$

and $T^{--} = P^-T\,|\,\mathcal{H}^-$ is a self-adjoint operator in $(\mathcal{H}^-, -[\cdot\,,\cdot])$, and the identity operator $I^- = I\,|\,\mathcal{H}^-$ (and similarly $I^+ = I\,|\,\mathcal{H}^+$).

Let $y^- = -P^-T^cx^+$, $x^+ \in \mathcal{H}^+ \cap \mathcal{D}_{T^c}$. By Lemma 2.4 d)

$$\|y^-\|^2 = \langle T_s y^-\,,\, x^+\rangle \le m\|y^-\|\,\|x^+\|$$

by Cauchy–Schwarz. Thus $\|y^-\| \le m\|x^+\|$ and

$$\|\tilde{K}_\lambda x^+\| \le \rho\|x^+\|\,, \quad \rho = \frac{m}{|\operatorname{Im}\lambda|}\,.$$

Because $\|K_\lambda\| \le \|\tilde{K}_\lambda\| \le \rho$, $\mathcal{N}_\lambda$ is uniformly positive if $\rho < 1$.

Consider some other canonical decomposition $\mathcal{H} = \mathcal{H}_1^-[\oplus]\mathcal{H}_1^+$ with the associated canonical symmetry $J_1$. Let $(x, \lambda J_1 x) \in T^c$, $x = x^- + x^+$, $x^\pm \in \mathcal{H}^\pm$, and $\lambda \in \mathbb{C}_*$. Since ([6, Theorem 1.8.17]) $J_1 x = Gx^- + Hx^+$ where

$$\begin{aligned} G =& -(I^+ + K)(I^+ - K^*K)^{-1}K^* - (I^- + K^*)(I^- - KK^*)^{-1}\,, \\ H =& (I^+ + K)(I^+ - K^*K)^{-1} + (I^- + K^*)(I^- - KK^*)^{-1}K \end{aligned}$$

where $K$ is the angular operator for $\mathcal{H}_1^+$ with respect to $\mathcal{H}^+$ (and $K^*$ is its Hilbert space adjoint, *i.e.* the angular operator for $\mathcal{H}_1^-$ with respect to $\mathcal{H}^-$), $(x, \lambda P^- J_1 x) \in P^- T^c$ equivalently reads thus

$$\begin{aligned} 0 =& (T^{--} - \lambda P^- G)x^- \\ & + ((P^- T^c \,|\, \mathcal{H}^+) - \lambda K + \lambda P^- GK)x^+ \end{aligned} \tag{2.2}$$

where one uses $P^- H = K - P^- GK$.

By extension by continuity $\tilde{K}_\lambda$ admits an extension on $\overline{\mathcal{H}^+ \cap \mathcal{D}_{T^c}}$, which we denote by the same $\tilde{K}_\lambda$. For $\rho < 1$, we define a strict contraction $K$ on $\mathcal{H}^+$ as the trivial extension ([6, Remark 1.8.8]) of $\tilde{K}_{-\lambda}$: $Kx^+ = \tilde{K}_{-\lambda}x^+$ if $x^+ \in \overline{\mathcal{H}^+ \cap \mathcal{D}_{T^c}}$; $Kx^+ = 0$ if $x^+ \in \mathcal{H}^+ \ominus \overline{\mathcal{H}^+ \cap \mathcal{D}_{T^c}}$ $(= \operatorname{Ind} T)$. With this $K$ (2.2) reads

$$0 = (T^{--} - \lambda P^- G)(x^- - Kx^+)$$

and it remains to verify that $T^{--} - \lambda P^- G$ has trivial kernel if $\lambda \in \Gamma_T$.

Using the definition of $G$, by elementary manipulations

$$\operatorname{Ker}(T^{--} - \lambda P^- G) = \operatorname{Ker}((I^- - KK^*)(T^{--} - \lambda I^-) + 2\lambda I^-)\,.$$

But the operator norm

$$\|(I^- - KK^*)(T^{--} - \lambda I^-)\| \leq (1 + \rho^2)|\|T^{--}\| \pm \lambda|$$

if $\operatorname{Re}\lambda \geq 0 (< 0)$, by the spectral theorem. □

*Remark* 2.6. If $\operatorname{Ker}_\lambda J_1 T^c \subseteq \mathcal{H}_1^+$ as in Lemma 2.5, to some other $\lambda_0 \in \Gamma_T$ there corresponds some other canonical decomposition $\mathcal{H} = \mathcal{H}_0^-[\oplus]\mathcal{H}_0^+$ with canonical symmetry $J_0$ such that $\operatorname{Ker}_{\lambda_0} J_0 T^c \subseteq \mathcal{H}_0^+$. Thus, for example $\operatorname{Ker}_{\bar{\lambda}} J_1 T^c \subseteq \tilde{V}(\mathcal{H}_1^+)$, $\tilde{V} = \left(\begin{smallmatrix} V^- & 0 \\ 0 & I^+ \end{smallmatrix}\right)$ with respect to $\mathcal{H} = \mathcal{H}^-[\oplus]\mathcal{H}^+$, $V^- = (T^{--} - \lambda I^-)(T^{--} - \bar{\lambda} I^-)^{-1}$ (Cayley transform).

2.5 Let $T$ be a relation in a $J$-space $\mathcal{H} = \mathcal{H}^-[\oplus]\mathcal{H}^+$. One verifies

$$\begin{aligned} P^-(JT - \lambda I)^{-1} \,|\, \mathcal{H}^- =& -(\mathcal{T}^-(\lambda) + \lambda I^-)^{-1}\,, \\ P^-(T - \lambda I)^{-1} \,|\, \mathcal{H}^- =& -(\mathcal{T}^-(\lambda) - \lambda I^-)^{-1} \end{aligned}$$

for all $\lambda \in \mathbb{C}$ $(I^- = I \,|\, \mathcal{H}^-)$, where the relation $\mathcal{T}^-(\lambda)$ in $\mathcal{H}^-$ is defined by

$$\mathcal{T}^-(\lambda) = \{(P^-x, P^-y) \,|\, (x, y) \in T\,;\; y - \lambda x \in \mathcal{H}^-\}$$

and is characterized in the next lemma.

**Lemma 2.7.** a) *If $\mathcal{H}^- \subseteq \mathcal{D}_{T^c}$ then $\mathcal{T}^-(\lambda)$ $(\lambda \in \mathbb{C} \setminus \sigma_p(P^+T \,|\, \mathcal{H}^+))$ is an operator.*

b) *If $T = T^c$ then $\mathcal{T}^-(\lambda)$ $(\lambda \in \mathbb{C}_*)$ is a Nevanlinna family in a Hilbert space $(\mathcal{H}^-, -[\cdot\,,\,\cdot])$.*

c) *Let $T_0 \supseteq T$ in $\mathcal{H}$ and let $\mathcal{T}_0^-$ be defined similar to $\mathcal{T}^-$, but with $T$ replaced by $T_0$. For $\lambda \in \mathbb{C} \setminus (\sigma_p(JT_0) \cup \sigma_p(P^+T_0 \,|\, \mathcal{H}^+))$, $\mathcal{T}^-(\lambda) \subseteq \mathcal{T}_0^-(\lambda)$ becomes the equality if $\mathcal{H}^- \cap \operatorname{Ind} T = \mathcal{H}^- \cap \operatorname{Ind} T_0$ and $\mathcal{H}^- \cap \mathcal{R}_{JT-\lambda I} = \mathcal{H}^- \cap \mathcal{R}_{JT_0-\lambda I}$.*

d) *If $\lambda \in \mathbb{C} \setminus (\{0\} \cup \sigma_p(JT))$ then $\lambda \in \sigma_p(T)$ iff $0 \in \sigma_p(\mathcal{T}^-(\lambda) - \lambda I^-)$.*

To recall (*e.g.* [11]), a Nevanlinna family (function) $\mathcal{T}(\lambda)$ is a relation (operator) in a Hilbert space, $(\mathcal{G}, \langle\cdot\,,\,\cdot\rangle_{\mathcal{G}})$ say, such that: a) $\operatorname{Im}\langle x\,,\,y\rangle_{\mathcal{G}} \geq 0$ $(\leq 0)$ for all $(x, y) \in \mathcal{T}(\lambda)$ and $\lambda \in \mathbb{C}^+$ $(\mathbb{C}^-)$; b) $\mathcal{T}(\lambda)^* = \mathcal{T}(\bar{\lambda})$ for all $\lambda \in \mathbb{C}_*$; c) $\lambda \mapsto (\mathcal{T}(\lambda) + \lambda_0 I)^{-1}$ is a holomorphic function on $\mathbb{C}^+$ $(\mathbb{C}^-)$ with values in $\mathcal{B}(\mathcal{G})$, the set of bounded everywhere defined operators in $\mathcal{G}$, for some and then for all $\lambda_0 \in \mathbb{C}^+$ $(\mathbb{C}^-)$.

*Proof.* a) By definition $x^- \in \operatorname{Ind}\mathcal{T}^-(\lambda) = \mathcal{H}^- \cap (JT - \lambda I)(\mathcal{H}^+)$ implies $x^- \in \mathcal{H}^-$ and $(\exists y^+ \in \operatorname{Ker}_\lambda P^+T \,|\, \mathcal{H}^+)$ $(y^+, \lambda y^+ + x^-) \in JT$, so $x^- \in \mathcal{H}^- \cap \operatorname{Ind} JT = \{0\}$ if $\lambda \notin \sigma_p(P^+T \,|\, \mathcal{H}^+)$.

b) $\operatorname{Im}\langle P^-x\,,\,P^-y\rangle = \operatorname{Im}\lambda \|P^+x\|^2$ (Section 2.3) is clear for $(x, y) \in T$, $y - \lambda x \in \mathcal{H}^-$, by using $\operatorname{Im}[x\,,\,y] = 0$.

By arguing as in Lemma 2.4 a) $P^-(JT-\lambda I)^{-1} \,|\, \mathcal{H}^-$ $(\lambda \in \mathbb{C}_*)$ is a closed operator in $(\mathcal{H}^-, -[\cdot\,,\,\cdot])$; since the adjoint $P^-(JT - \bar{\lambda} I)^{-1} \,|\, \mathcal{H}^-$ equals both $-(\mathcal{T}^-(\lambda)^* + \bar{\lambda} I^-)^{-1}$ and $-(\mathcal{T}^-(\bar{\lambda}) + \bar{\lambda} I^-)^{-1}$, so $\mathcal{T}^-(\lambda)^* = \mathcal{T}^-(\bar{\lambda})$.

Finally, $(\mathcal{T}^-(\lambda) + \lambda I^-)^{-1} \in \mathcal{B}(\mathcal{H}^-)$ for $\lambda \in \mathbb{C}_* \subseteq \rho(JT)$.

c) Since $\mathcal{T}^-(\lambda) \subseteq \mathcal{T}_0^-(\lambda)$ $(\lambda \in \mathbb{C})$, the equality holds if the domains and multivalued parts of $\mathcal{T}^-(\lambda)$ and $\mathcal{T}_0^-(\lambda)$ coincide; as already shown in a) $\operatorname{Ind}\mathcal{T}^-(\lambda) = \mathcal{H}^- \cap \operatorname{Ind} T$, and similarly $\operatorname{Ind}\mathcal{T}_0^-(\lambda) = \mathcal{H}^- \cap \operatorname{Ind} T_0$, if $\lambda \notin \sigma_p(P^+T_0 \,|\, \mathcal{H}^+)$.

Assume that $\mathcal{H}^- \cap \mathcal{R}_{JT-\lambda I} = \mathcal{H}^- \cap \mathcal{R}_{JT_0-\lambda I}$ for some $\lambda \in \mathbb{C}$. Applying $(JT_0-\lambda I)^{-1}$ on both sides gives $\operatorname{Ker}_\lambda JT_0 + (JT-\lambda I)^{-1}(\mathcal{H}^-) = (JT_0-\lambda I)^{-1}(\mathcal{H}^-)$. By noting that $\mathcal{D}_{\mathcal{T}^-(\lambda)} = P^-(JT-\lambda I)^{-1}(\mathcal{H}^-)$, and similarly for $\mathcal{D}_{\mathcal{T}_0^-(\lambda)}$, the assertion follows.

d) This follows from

$\operatorname{Ker}_\lambda \mathcal{T}^-(\lambda)) = P^-(\operatorname{Ker}_\lambda T)$ $(\lambda \in \mathbb{C})$;

$\operatorname{Ker}_\lambda JT + \operatorname{Ker}_\lambda T = (JT - \lambda I)^{-1}(\operatorname{Ker}_\lambda \mathcal{T}^-(\lambda))$ $(\lambda \in \mathbb{C} \setminus \{0\})$.

The first formula is clear from the definition of $\mathcal{T}^-$, so we only prove the second one. Let $u \in (JT - \lambda I)^{-1} P^-(\operatorname{Ker}_\lambda T)$; hence $(u, x^-) \in JT - \lambda I$, $x^- = P^-x$, $(x, \lambda x) \in T$. Since $(2\lambda u, 2\lambda x^-) \in JT - \lambda I$ on the one hand and $(x, 2\lambda x^-) \in JT - \lambda I$ on the other hand, $2\lambda u - x \in \operatorname{Ker}_\lambda JT$.

Conversely, $u = x + y$ and $(x, J\lambda x) \in JT$ and $(y, \lambda y) \in JT$ implies $(u, -2\lambda x^-) \in JT - \lambda I$, $x^- = P^- x$. But $-2\lambda x \in \operatorname{Ker}_\lambda T$. $\square$

2.6 Proof of Theorem 1.3 a): Fix an arbitrary $\lambda \in \mathbb{C}^- \cap \Gamma_T$. By Lemma 2.5 one finds a canonical decomposition $\mathscr{H} = \mathscr{H}_1^- [\oplus] \mathscr{H}_1^+$ with canonical symmetry $J_1$ such that $\mathscr{R}_{J_1 T - \lambda I} \supseteq \mathscr{H}_1^-$. Because $\operatorname{Ind} T_0 \subseteq \mathscr{H}^+$, $\mathscr{H}_1^- \cap \operatorname{Ind} T_0 = \{x^- + K^* x^- \mid x^- \in \mathscr{H}^-\} \cap \operatorname{Ind} T_0 = \{0\}$ (where $K^* = P^+(P^- \mid \mathscr{H}_1^-)^{-1}$), so $\lambda \in \sigma_p(T_0)$ iff $\lambda \in \sigma_p(T)$ by Lemma 2.7 c), d). From here, replacing $T$ ($T_0$) by $T \mid \mathscr{H}^-$ ($T$), $\Gamma_T \cap \sigma_p(T) = \Gamma_T \cap \sigma_p(T \mid \mathscr{H}^-) = \emptyset$ by Lemma 2.4 a). Thus $\mathbb{C}^- \cap \sigma_p(T_0)$ is contained in $\mathbb{C}^- \setminus \Gamma_T (\subseteq C_T)$. $\square$

The proof shows: the only place where one uses dissipativeness is when applying Lemma 2.7 c), in order to get rid of $\sigma_p(J_1 T_0) \cup \sigma_p(P_1^+ T_0 \mid \mathscr{H}_1^+)$. Possible generalizations, which we leave out, are now evident.

## 3. Decomposition of point spectrum of the componentwise sum of two relations

3.1 Let $T$ and $N$ be relations in a linear space. Their componentwise sum (aka linear span) is denoted $T \hat{+} N$. We are interested in $\operatorname{Ker}_\lambda(T \hat{+} N)$ for some $\lambda \in \mathbb{C}$. Formally

**Proposition 3.1.** (*cf.* [12, Lemma 1.7.2])

$$\begin{aligned}\operatorname{Ker}_\lambda(T \hat{+} N) &= \mathscr{R}_{(T-\lambda I)^{-1} - (N-\lambda I)^{-1}} \\ &= \mathscr{R}_{(T-\lambda I)^{-1}(\lambda I - N) + I} \\ &= \mathscr{R}_{(T-\lambda I)^{-1}(N-\lambda I) - I} \,.\end{aligned}$$

*Proof.* $\operatorname{Ker}_\lambda(T \hat{+} N)$ consists of those $x - y$ such that $(\exists u \in \mathscr{R}_{T-\lambda I} \cap \mathscr{R}_{N-\lambda I})$ $(x, u) \in T - \lambda I$ and $(y, u) \in N - \lambda I$. $\square$

3.2 The aim, however, is to extract $\operatorname{Ker}_\lambda T + \operatorname{Ker}_\lambda N$ from $\operatorname{Ker}_\lambda(T \hat{+} N)$. The results in this section, therefore, should be seen as a further investigation of some of the ideas presented in [30].

**Proposition 3.2.** *In order that*

$$\operatorname{Ker}_\lambda(T \hat{+} N) = \operatorname{Ker}_\lambda T + \operatorname{Ker}_\lambda N \tag{3.1}$$

*should hold for some* $\lambda \in \mathbb{C}$ *it is necessary and sufficient that*

$$\begin{aligned} &(N - \lambda I)^{-1}(\mathscr{R}_{T-\lambda I}) = \mathscr{D}_{T \cap N} + \operatorname{Ker}_\lambda N \\ \textit{and} \quad &\mathscr{R}_{T-\lambda I} \cap \operatorname{Ind} N = ((T \cap N) - \lambda I)(\operatorname{Ker}_\lambda N) \,. \end{aligned}$$

*Proof.* Fix a complex number $\lambda$ and define a relation

$$Z = \{(x, x') \in T \mid x' - \lambda x \in \mathscr{R}_{N-\lambda I}\} \,.$$

Since

$$Z \supseteq T \cap N \,, \quad Z \supseteq \lambda I \cap T$$

by linearity

$$Z \supseteq Z_0 = (T \cap N) \mathbin{\hat{+}} (\lambda I \cap T)\,.$$

(3.1) holds iff $Z = Z_0$:

Because $\mathrm{Ker}_\lambda(T \mathbin{\hat{+}} N)$ consists of those $x-y$ such that $(\exists x')\ (x,x') \in Z$ and $(y, x' + \lambda(y-x)) \in N$, if $Z = Z_0$ then $(x,x') = (u+v, u'+\lambda v)$, $(u,u') \in T \cap N$, $v \in \mathrm{Ker}_\lambda T$, so $x - y = v + (u-y)$, $u - y \in \mathrm{Ker}_\lambda N$.

Conversely, if (3.1), consider $x - y \in \mathrm{Ker}_\lambda(T \mathbin{\hat{+}} N)$ as described above; because also $x - y = x_\lambda - y_\lambda$ for some $x_\lambda \in \mathrm{Ker}_\lambda T$ and $y_\lambda \in \mathrm{Ker}_\lambda N$, one gets $(x - x_\lambda, x' - \lambda x) \in (T \cap N) - \lambda I$ *i.e.* $(x,x') \in Z_0$.

There is the equality in $Z \supseteq Z_0$ iff

$$\begin{aligned} &(T-\lambda I)^{-1}(\mathscr{R}_{N-\lambda I}) = \mathscr{D}_{T\cap N} + \mathrm{Ker}_\lambda T \\ \text{and}\quad &\mathscr{R}_{N-\lambda I} \cap \mathrm{Ind}\, T = ((T \cap N) - \lambda I)(\mathrm{Ker}_\lambda T) \end{aligned}$$

(*i.e.* their domains and multivalued parts coincide). Relabeling $T$ $(N)$ by $N$ $(T)$, the assertion follows. □

We note without proof that

$$((T \cap N) - \lambda I)(\mathrm{Ker}_\lambda N) = (\mathscr{R}_{(T\cap N)-\lambda I} \cap \mathrm{Ker}_\lambda N) + \mathrm{Ind}(T \cap N)\,.$$

3.3 Let $T$, $N$ be as in the previous paragraph. Let $O = O(T,N)$ denote the set of those $\lambda \in \mathbb{C}$ such that

$$\mathscr{R}_{T-\lambda I} \cap \mathscr{R}_{N-\lambda I} = \mathrm{Ind}(T \cap N)\,. \tag{3.2}$$

Let $O^c = O^c(T,N) = \mathbb{C} \setminus O$.

3.4 As a corollary of Proposition 3.2, in order that a number $\lambda$ should be in $O$ it is necessary that $\mathscr{D}_{T\cap N} = \mathrm{Ker}_\lambda(T \cap N)$. And conversely, if the latter holds, (3.1) and (3.2) are equivalent. Therefore

**Corollary 3.3.** *If $T \cap N$ has trivial domain, then*

$$\mathbb{C} \setminus \sigma_p(T \mathbin{\hat{+}} N) = O \setminus \sigma_p(\begin{pmatrix} T & 0 \\ 0 & N \end{pmatrix})\,.$$

*In particular*

$$\sigma_p(T \mathbin{\hat{+}} N) = \sigma_p(\begin{pmatrix} T & 0 \\ 0 & N \end{pmatrix}) \cup O^c\,.$$

*Remarks* 3.4. a) The situation changes drastically if $\mathscr{D}_{T\cap N}$ is nontrivial: Since every $\lambda \in \mathbb{C} \setminus \sigma_p(T \mathbin{\hat{+}} N)$ satisfies (3.1), $O \subseteq \sigma_p(T \mathbin{\hat{+}} N)$ by Proposition 3.2. If also there is no $\lambda$ with $\mathscr{D}_{T\cap N} = \mathrm{Ker}_\lambda(T \cap N)$, $O$ is empty. Take *e.g.* $N \subseteq T$, $\mathscr{D}_N \neq \{0\}$, $\sigma_p(N) = \emptyset$.

b) $E = \left(\begin{smallmatrix} T & 0 \\ 0 & N \end{smallmatrix}\right)$ and $T_0 = T \mathbin{\hat{+}} N$ are related by $T_0 = VEV^{-1}$, $V\colon \binom{x}{y} \mapsto x + y$. From here $\mathrm{Ker}_\lambda T_0 = V(E - \lambda I)^{-1}(\mathrm{Ker}\, V)$. Thus $\mathrm{Ker}_\lambda T_0 = V(\mathrm{Ker}_\lambda E)$ *i.e.* (3.1) holds if $\mathscr{R}_{E-\lambda I} \cap \mathrm{Ker}\, V \subseteq \mathrm{Ind}\, E$ or, what is the same, $\lambda \in O$. For block relations we refer to [26].

*Example* 3.5. By Proposition 3.2, for $\lambda$ to be in $O(\tilde{T},\tilde{N})$ it is sufficient—and if $\mathcal{D}_T\cap\mathcal{D}_N=\{0\}$ then also necessary—that

$$(\mathcal{D}_T+\operatorname{Ker}_\lambda(T\,\hat{+}\,N))\cap\mathcal{D}_N=\operatorname{Ker}_\lambda\tilde{N}$$

where $\tilde{N}=N\,\hat{+}\,(\{0\}\times\mathcal{R}_{T-N})$ and similarly for $\tilde{T}$. For example, let $T$ be a closed densely defined symmetric operator in a Hilbert space. Consider a closed operator $T_0\supseteq T$: $T_0=T\,\hat{\oplus}\,N$, $\mathcal{D}_N=\operatorname{Ker}_{-1}T^*T_0$. Then

$$\begin{aligned}O=&(\mathbb{C}\setminus\sigma_p(T_0))\sqcup(\{-\mathrm{i},\mathrm{i}\}\cap\sigma_p(T_0))\\&\sqcup\{\lambda\in\sigma_p(T)\mid\operatorname{Ker}_\lambda T_0=\operatorname{Ker}_\lambda T\}\,.\end{aligned}$$

This $O$ also follows from the next lemma.

**Lemma 3.6.** *Let* $T\cap N=\{0\}$, $\lambda\in\mathbb{C}$. *Then* a)–c) *are equivalent:*

a) $\lambda\in O$;

b) $(T\,\hat{+}\,[\lambda I\cap(T\,\hat{+}\,N)])\cap N\subseteq\lambda I$;

c) $(\mathcal{D}_T+\operatorname{Ker}_\lambda(T\,\hat{+}\,N))\cap\mathcal{D}_N=\operatorname{Ker}_\lambda N$, $\mathcal{R}_{T-\lambda I}\cap\operatorname{Ind}N=\{0\}$.

*Therefore* $\sigma_p(T\,\hat{+}\,N)=\mathbb{C}$ *if* $\mathcal{D}_T\cap\mathcal{D}_N$ *is nontrivial.*

*Proof.* Let $T_0=T\,\hat{+}\,N$.

a) $\Rightarrow$ b) Let $(x,y)\in(T\,\hat{+}\,(\lambda I\cap T_0))\cap N$, $\lambda\in O$, then $(\exists h\in\operatorname{Ker}_\lambda T_0)$ $(x-h,y-\lambda h)\in T$, *i.e.* $y-\lambda x\in\mathcal{R}_{T-\lambda I}\cap\mathcal{R}_{N-\lambda I}=\{0\}$.

a) $\Leftarrow$ b) Let $(x,y)\in N$, $(u,v)\in T$, $y-\lambda x=v-\lambda u$, $\lambda\in\mathbb{C}$. Then $(x-u,\lambda(x-u))\in T_0$ and therefore $(x,y)=(u,v)+(x-u,\lambda(x-u))\in(T\,\hat{+}\,(\lambda I\cap T_0))\cap N$.

b) $\Leftrightarrow$ c) By [48, Proposition 1], b) holds iff $(T\,\hat{+}\,(\lambda I\cap T_0))\cap N$ is an operator (*i.e.* $\mathcal{R}_{T-\lambda I}\cap\operatorname{Ind}N=\{0\}$) with domain $(\mathcal{D}_T+\operatorname{Ker}_\lambda T_0)\cap\mathcal{D}_N=\operatorname{Ker}_\lambda N$.

To prove that last statement consider $\lambda\in O\setminus\sigma_p(N)$. Then by c) either necessarily $\mathcal{D}_T\cap\mathcal{D}_N=\{0\}$ or else $O\subseteq\sigma_p(N)$, and in the latter case then necessarily $\sigma_p(T_0)\supseteq O^c\supseteq\mathbb{C}\setminus\sigma_p(N)\supseteq\mathbb{C}\setminus\sigma_p(T_0)$ *i.e.* $\sigma_p(T_0)=\mathbb{C}$. □

## 4. Deficiency subspace and related topics

4.1 The adjoint of a symmetric relation $T$ in a $J$-space $\mathcal{H}$ decomposes as

$$T^c=\bar{T}\;[\hat{\oplus}]\;\Sigma$$

where the closed relation $\Sigma=T^c\cap T^\perp$ in $\mathcal{H}$ explicitly reads

$$\begin{aligned}\Sigma=&\{(x,x')\in T^c\mid Jx'+\mathrm{i}x\in\operatorname{Ker}_{\mathrm{i}}JT^c\}\\=&\{(x,x')\in T^c\mid Jx'-\mathrm{i}x\in\operatorname{Ker}_{-\mathrm{i}}JT^c\}\\=&(\mathrm{i}J\mid\operatorname{Ker}_{\mathrm{i}}JT^c)\;[\hat{\oplus}]\;(-\mathrm{i}J\mid\operatorname{Ker}_{-\mathrm{i}}JT^c)\,.\end{aligned}\tag{4.1}$$

The linear span $\Sigma$ of deficiency subspaces $\pm\mathrm{i}J\mid\operatorname{Ker}_{\pm\mathrm{i}}JT^c$ ([6, Proposition 1.4.20]) is itself called here the deficiency subspace.

Whenever we discuss $\Sigma$, by (4.1) it is no loss of generality to assume that a symmetric $T$ is closed.

*Remark* 4.1. $(\Sigma, [\cdot, \cdot]_\Gamma)$ is a Krein subspace of $(\mathcal{H}_\Gamma, [\cdot, \cdot]_\Gamma)$, and the last equality in (4.1) is its canonical decomposition.

4.2 $\Sigma^c = T \,\widehat{\oplus}\, J_\Gamma(T)$ is the $J$-adjoint of the deficiency subspace $\Sigma$ of a closed $J$-symmetric relation $T$. A few elementary characteristics:
a) $\mathcal{D}_\Sigma = \mathrm{Ker}_{-1}(JT^c)^2 = \mathrm{Ker}_{\mathrm{i}}\, JT^c + \mathrm{Ker}_{-\mathrm{i}}\, JT^c = J(\mathcal{R}_\Sigma)$.
b) $\mathcal{D}_\Sigma^\perp = \mathcal{R}_{(JT)^2+I} = \mathcal{R}_{JT+\mathrm{i}I} \cap \mathcal{R}_{JT-\mathrm{i}I}$.
c) $\mathcal{D}_{\Sigma^c} = \mathcal{D}_T + J(\mathcal{R}_T) = J(\mathcal{R}_{\Sigma^c})$.

By c) equivalent are: $\Sigma$ is an operator; $JT \in (P)$; $\Sigma$ is injective. If $\bar{\mathcal{D}}_T$ is regular ([6, Proposition 1.7.15]), then one adds an additional equivalence $T \in (P)$. The condition $(P)$ is characterized by the next lemma.

**Lemma 4.2.** ([30, Lemma A.5]) *Let $T$ be a relation in a Krein space. If* $\mathrm{Ker}_\lambda T^c \cap \mathrm{Ker}_{\lambda_0} T^c$ *is trivial for some $\lambda$, $\lambda_0 \in \sigma_p(T^c)$, $\lambda \neq \lambda_0$, then it is trivial for all $\lambda$, $\lambda_0 \in \sigma_p(T^c)$, $\lambda \neq \lambda_0$. In order that this should happen it is necessary and sufficient that $T \in (P)$.*

In particular $\{\mathrm{Ker}_{\mathrm{i}}\, JT^c, \mathrm{Ker}_{-\mathrm{i}}\, JT^c\}$ is a linearly independent system for a $J$-symmetric $T \in (P)$ $(T \neq T^c)$ with $\bar{\mathcal{D}}_T$ regular.

*Remark* 4.3. Among other things, a relation $T$ such that $\mathcal{D}_T + \mathcal{R}_T = \mathcal{H}$ in a Hilbert space is examined *e.g.* in [43] (but note that a surjective symmetric operator is self-adjoint, *e.g.* [48, Proposition 2], [44, Corollary 3.12]).

4.3 Let $T \neq T^c$ be a closed symmetric relation in a $J$-space $(\mathcal{H})$. Although each separately $\{\mathcal{D}_T, \mathrm{Ker}_{\mathrm{i}}\, JT^c\}$, $\{\mathcal{D}_T, \mathrm{Ker}_{-\mathrm{i}}\, JT^c\}$ are two linearly independent systems (*e.g.* [22, Eq. (1.2)]), that is not necessarily the case for the system $\{\mathcal{D}_T, \mathrm{Ker}_{\mathrm{i}}\, JT^c, \mathrm{Ker}_{-\mathrm{i}}\, JT^c\}$: By Section 4.2 a) $\mathcal{D}_T \cap \mathcal{D}_\Sigma = \{0\}$ iff $-1 \notin \sigma_p(T^*T_c)$, where $T_c = T^c \,|\, \mathcal{D}_T$ in $\mathcal{H}$, see (2.1), and $T^* = JT^cJ$ is the Hilbert adjoint. The next lemma uses another viewpoint.

**Lemma 4.4.** (*cf.* [32, Theorem 8]) *Let $T$ be closed $J$-symmetric. Then*

$$\mathcal{D}_T \cap \mathcal{D}_\Sigma = (P_{\mathrm{i}} - P_{-\mathrm{i}})(\mathrm{Ind}\, JT^c)$$

*where $P_{\pm\mathrm{i}}$ are orthogonal projections in $(\mathcal{H}, [\cdot, J\cdot])$ onto* $\mathrm{Ker}_{\pm\mathrm{i}}\, JT^c$.

*Proof.* Consider the equation $x + x_{\mathrm{i}} + x_{-\mathrm{i}} = 0$, $x \in \mathcal{D}_T$, $x_{\pm\mathrm{i}} \in \mathcal{N}_{\pm\mathrm{i}} = \mathrm{Ker}_{\pm\mathrm{i}}\, JT^c$. Using $(-x_{\mathrm{i}} - x_{-\mathrm{i}}, x') \in JT$ $(\subseteq JT^c)$, some $x' \in \mathcal{H}$, yields $x' + \mathrm{i}x_{\mathrm{i}} - \mathrm{i}x_{-\mathrm{i}} = y \in \mathrm{Ind}\, JT^c$, and then

$$\begin{aligned} \mathrm{i}x_{\mathrm{i}} &= P_{\mathrm{i}}(y - x' + \mathrm{i}x_{-\mathrm{i}}) = P_{\mathrm{i}}(y - x' + \mathrm{i}(-x - x_{\mathrm{i}})) \\ &= P_{\mathrm{i}}y - P_{\mathrm{i}}(x' + \mathrm{i}x) - \mathrm{i}x_{\mathrm{i}} \end{aligned}$$

*i.e.* $2\mathrm{i}x_{\mathrm{i}} = P_{\mathrm{i}}y$ since $x' + \mathrm{i}x \in \mathcal{N}_{\mathrm{i}}^\perp (= \mathcal{R}_{JT+\mathrm{i}I})$. Similarly $2\mathrm{i}x_{-\mathrm{i}} = -P_{-\mathrm{i}}y$.

Conversely, consider $x = -x_{\mathrm{i}} - x_{-\mathrm{i}}$, $x_{\pm\mathrm{i}} \in P_{\pm\mathrm{i}}(\operatorname{Ind} JT^c)$; *i.e.* $(\exists y \in \operatorname{Ind} JT^c)$ $2\mathrm{i}x_{\pm\mathrm{i}} = \pm P_{\pm\mathrm{i}}y$. Letting $x' = y - \mathrm{i}x - 2\mathrm{i}x_{\mathrm{i}}$ yields

$$x' + \mathrm{i}x = y - 2\mathrm{i}x_{\mathrm{i}} = (I - P_{\mathrm{i}})y \in \mathscr{R}_{JT+\mathrm{i}I}\,,$$
$$x' - \mathrm{i}x = y + 2\mathrm{i}x_{-\mathrm{i}} = (I - P_{-\mathrm{i}})y \in \mathscr{R}_{JT-\mathrm{i}I}$$

*i.e.* $(\exists(u,u') \in JT)$ $(\exists(v,v') \in JT)$

$$x' + \mathrm{i}x = u' + \mathrm{i}u\,, \quad x' - \mathrm{i}x = v' - \mathrm{i}v\,.$$

Solving with respect to $(x,x')$ yields

$$(x,x') = (m,m') + (n,n')$$

for some $(m,m') \in JT$, $(n,n') \in -(JT)^{-1}(= (JT^c)^{\perp})$.

Now $(n,n') = (0,0)$, because $(x,x') \in JT^c$, $(m,m') \in JT^c$.

As a result $x \in \mathscr{D}_T$. □

**Proposition 4.5.** *Let $T$ be closed $J$-symmetric, $\bar{\mathscr{D}}_T$ regular. Then*

a) *$\Sigma$ is an operator iff $T \in (P)$.*

b) *If $\Sigma$ is an operator then $\mathscr{D}_T \cap \mathscr{D}_\Sigma$ is trivial iff $T \in (P')$.*
*More generally, if $\Sigma$ is an operator, let $T_0 \supseteq T$ be a closed extension in $\mathscr{H}$, and $N = T_0 \cap \Sigma$. Then $\mathscr{D}_T \cap \mathscr{D}_N = (\Sigma - T_s)^{-1}(\operatorname{Ind} T_0)$; hence the lineals $\mathscr{D}_T \cap \mathscr{D}_N$ and $\mathscr{R}_{\Sigma - T_s} \cap \operatorname{Ind} T_0$ are either both trivial or both nontrivial, and if $T$ is an operator, $\mathscr{D}_T \cap \mathscr{D}_N$ is trivial iff $T_0 \cap T^c$ is an operator.*

*Proof.* It remains to verify b).

If $\Sigma$ is an operator, *i.e.* $JT \in (P)$, then $\operatorname{Ker}_{\mathrm{i}} JT^c$, $\operatorname{Ker}_{-\mathrm{i}} JT^c$ are linearly independent (Lemma 4.2), so that $\mathscr{D}_T \cap \mathscr{D}_\Sigma = \{0\}$ iff $\operatorname{Ind} JT^c \subseteq \mathscr{R}_{(JT)^2+I}$ (Section 4.2 b) and Lemma 4.4) iff $\operatorname{Ind} T^c = \operatorname{Ind} T$ (*e.g.* [25, Eq. (2.4)]) iff $JT \in (P')$ iff $T \in (P')$.

As for the remaining part of b), because $N \subseteq \Sigma$ are operators, by definition $N = \{(x, \Sigma x) \in T_0 \mid x \in \mathscr{D}_{T_0} \cap \mathscr{D}_\Sigma\}$, so using the canonical form of $T_0$, with $(T_0)_s$ its operator part, $x \in \mathscr{D}_T \cap \mathscr{D}_N$ iff $x \in \mathscr{D}_H$ and $(\Sigma - H)x \in \operatorname{Ind} T_0$, where $H = (T_0)_s \mid \mathscr{D}_T \cap \mathscr{D}_\Sigma$. Because $H \subseteq T_s \hat{+} (\{0\} \times \operatorname{Ind} T_0)$, $(\Sigma - T_s)x \in \operatorname{Ind} T_0$, and this shows $\mathscr{D}_T \cap \mathscr{D}_N \subseteq (\Sigma - T_s)^{-1}(\operatorname{Ind} T_0)$.

For $\supseteq$, consider $x \in \mathscr{D}_T \cap \mathscr{D}_\Sigma$, $(\Sigma - T_s)x \in \operatorname{Ind} T_0$. Because $(x, T_s x) \in H \hat{+} (\{0\} \times \operatorname{Ind} T_0)$, so $(\Sigma - H)x \in \operatorname{Ind} T_0$ *i.e.* $x \in \mathscr{D}_T \cap \mathscr{D}_N$.

For the last two iff statements use that $(\Sigma - T_s)^{\pm 1}$ are operators and that $\mathscr{R}_{\Sigma - T} = \operatorname{Ind}(T \hat{+} \Sigma) = \operatorname{Ind} T^c$. □

*Remark* 4.6. Without using the condition $(P)$ and regularity of $\bar{\mathscr{D}}_T$, one easily verifies that the lineal $\mathscr{D}_T \cap \mathscr{D}_N$ can be given the form $\mathscr{R}_{V(T)\cap K(T_0)-I}$ in terms of relations in $\mathscr{H}$ defined by $V(T) = P_{\mathrm{i}}(P_{-\mathrm{i}} \mid \operatorname{Ind} JT^c)^{-1}$ (hence $\mathscr{D}_T \cap \mathscr{D}_\Sigma = \mathscr{R}_{V(T)-I}$ in Lemma 4.4) and $K(T_0) = \{(y - \mathrm{i}x, y + \mathrm{i}x) \mid (x,y) \in T_0\}$ (Cayley transform); if $T$ is an operator, so is $V(T)$, and which is studied in [32], and see also recent papers [3,4]. By adding conditions $T \in (P)$ and $\bar{\mathcal{D}}_T$ regular, therefore, the lineal $\mathcal{D}_T \cap \mathcal{D}_N$ is trivial iff such is the relation $V(T) \cap K(T_0)$. Just like Proposition 4.5 b), the criterion is of interest only if $T \notin (P')$.

By Lemma 3.6 and Proposition 4.5, if $T$ is a closed $J$-symmetric relation and $\bar{\mathcal{D}}_T$ is regular, then $\sigma_p(T^c) = \mathbb{C}$ if $T \in (P)$ but $T \notin (P')$. Moreover

**Corollary 4.7.** *Let $T$ be a closed $J$-symmetric operator, with regular $\bar{\mathcal{D}}_T$, and let $T_0$ be its closed extension in $\mathcal{H}$, with nonempty $\rho(T_0)$. If $T \in (P)$ then $T_0 \cap T^c$ is an operator.*

4.4 Proof of Theorem 1.3 b): By a) and Corollary 4.7. □

4.5 Conditions $(P)$, $(P')$ in Proposition 4.5 in terms of $T_s$ are reformulated as follows, where we use $\mathcal{L}_T = T^{-1}(\bar{\mathcal{D}}_T)$ for the domain of the range restriction $T \cap \bar{\mathcal{D}}_T^2$.

**Proposition 4.8.** *Let $T$ be closed $J$-symmetric, $\bar{\mathcal{D}}_T$ regular. Then*

a) *$T \in (P)$ iff $\overline{\mathcal{D}_T + \mathcal{R}_{T_s}} = \bar{\mathcal{D}}_{T^c}$ that is iff $T_s \in (P)$ in the Krein space $\mathcal{H}_c = (\bar{\mathcal{D}}_{T^c}, [\cdot, \cdot])$.*

b) *$T \in (P')$ iff $T_s \in (P')$ in $\mathcal{H}_c$.*

c) *$\mathcal{R}_{T_s} \subseteq \bar{\mathcal{D}}_T$ iff $\mathcal{L}_T = \mathcal{D}_T$; $\mathcal{L}_T = \mathcal{L}_{T_s}$.*

d) *If $T \in (P')$ then $T \in (P)$ with $\mathcal{L}_T = \mathcal{D}_T$.*

In the proposition, the operator $T_s$ in $\mathcal{H}_c$ is closed symmetric, with the adjoint given by the componentwise sum of the operator part of $T^c$ and the singular relation $\{0\} \times (\mathcal{D}_T [\oplus] \operatorname{Ind} T)^{[\perp]}$. If $T \in (P')$, the singular part vanishes, *i.e.* $T_s$ is a closed densely defined symmetric operator in the Krein space $\mathcal{H}_c = \bar{\mathcal{D}}_T$.

## 5. Class $(LP)$

5.1 Let $T$ be a relation in a Krein space $(\mathcal{H}, [\cdot, \cdot])$. By $\operatorname{Ext}(T)$, $\operatorname{Ext}^{\sim}(T)$, $\operatorname{Sym}(T)$ one means the set of canonical, closed canonical, closed symmetric extensions $T_0$: $T \subseteq T_0 \subseteq T^c$ respectively.

For short, $T$ is of class $(LP)$ (of class $(LP')$), if it is closed, symmetric, and has the properties $(L) \wedge (P)$ (resp. $(L) \wedge (P')$).

The spectrum (point spectrum, see Corollary 3.3, [46, Theorem 2.7]) of a closed symmetric relation $T \in (L)$, in particular $T \in (LP)$, coincides with that of $T_s$. If $T \in (LP')$, then also the (point) spectrum of $T_0 \in \operatorname{Ext}^{\sim}(T)$ coincides with the (point) spectrum of its operator part, the latter being a closed densely defined canonical extension of $T_s$ in $\mathcal{H}_c = \bar{\mathcal{D}}_T$.

5.2 Combined with Proposition 4.8 the next criterion is useful for establishing $T \in (LP)$, $T \notin (LP')$, if $\mathcal{L}_T$ is non-trivial.

**Proposition 5.1.** *Consider a closed symmetric relation $T\in(L)$ in a Krein space with canonical decomposition $\mathcal{H}=\mathcal{H}^-[\oplus]\mathcal{H}^+$; $\mathcal{H}^-\subseteq\mathcal{D}_T$. Then* a)–c) *are equivalent:*

a) $\mathcal{H}^-\subseteq\mathcal{L}_T$;

b) $T_s(\mathcal{H}^-)\subseteq\bar{\mathcal{D}}_T$;

c) $T_s(\mathcal{H}^-)\subseteq\overline{\mathcal{D}_T+\operatorname{Ind}T}$.

*If* a)–c) *hold then $\mathcal{L}_T=\mathcal{D}_T$ iff $T_s(\mathcal{H}^+)\subseteq\overline{\mathcal{D}_T+\operatorname{Ind}T}$.*

*Proof.* In b) $\Leftrightarrow$ c) it suffices to remark on b) $\Leftarrow$ c): $P^+(\mathcal{R}_{T_s})\subseteq\mathcal{H}_c^+$ where $\mathcal{H}_c^+=\mathcal{H}^+\cap\mathcal{H}_c$ (where $\mathcal{H}_c$ is as in Proposition 4.8).

a) $\Leftrightarrow$ b): The matrix representation of the operator $T_s$ in the Krein space $\mathcal{H}_c=\mathcal{H}^-[\oplus]\mathcal{H}_c^+$ is given by

$$T_s=\begin{pmatrix}T_s^{--} & T_s^{-+}\\ T_s^{+-} & T_s^{++}\end{pmatrix}$$

with the entries

$$\begin{aligned}T_s^{--}&=P^-T_s\,|\,\mathcal{H}^-=P^-T\,|\,\mathcal{H}^-\,,\\ T_s^{-+}&=P^-T_s\,|\,\mathcal{H}^+=P^-T\,|\,\mathcal{H}^+\,,\\ T_s^{+-}&=P_c^+T_s\,|\,\mathcal{H}^-=P^+T_s\,|\,\mathcal{H}^-=P_c^+T\,|\,\mathcal{H}^-\,,\\ T_s^{++}&=P_c^+T_s\,|\,\mathcal{H}^+=P^+T_s\,|\,\mathcal{H}^+=P_c^+T\,|\,\mathcal{H}^+\end{aligned}$$

where $P_c^+$ is the orthogonal projection in $\mathcal{H}$ onto $\mathcal{H}_c^+$.

Since $\mathcal{L}_T=\mathcal{L}_{T_s}$, $\mathcal{L}_T$ contains only those $x^-+x^+\in\mathcal{H}^-[\oplus](\mathcal{H}^+\cap\mathcal{D}_T)$ such that $T_s^{+-}x^-+T_s^{++}x^+\in\mathcal{H}^+\cap\bar{\mathcal{D}}_T$, and the assertion follows.

From the above also the last statement is evident. □

5.3 If $T\in(LP)$ in a Krein space $\mathcal{H}$ and $T_0\in\operatorname{Ext}^{\sim}(T)$, then the non-real point spectrum of $N=T_0\cap\Sigma$ is at most $\{\pm\mathrm{i}\}$.

**Proposition 5.2.** *The deficiency subspace $\Sigma$ of $T\in(LP)$ in a Krein space $\mathcal{H}$ is the operator with at most two non-real eigenvalues, $\pm\mathrm{i}$; $\operatorname{Ker}_{\pm\mathrm{i}}\Sigma=\mathcal{H}^+\cap\operatorname{Ker}_{\pm\mathrm{i}}JT^c$ in the canonical decomposition $\mathcal{H}=\mathcal{H}^-[\oplus]\mathcal{H}^+$, $\mathcal{H}^-\subseteq\mathcal{D}_T$, with canonical symmetry $J$.*

*If $T\in(LP')$ then $\Sigma$ has no non-real eigenvalues.*

*Proof.* Using (4.1) $\operatorname{Ker}_\lambda\Sigma$ ($\lambda\neq0,\pm\mathrm{i}$) consists of $x_{\mathrm{i}}-x_{-\mathrm{i}}$ where $x_{\pm\mathrm{i}}\in\operatorname{Ker}_{\pm\mathrm{i}}JT^c$, $(J-\mu I)x_{\mathrm{i}}=(J+\mu I)x_{-\mathrm{i}}$, $\mu=\mathrm{i}/\lambda\neq0,\pm1$. (If $\mu=0$ then $x_{\mathrm{i}}=x_{-\mathrm{i}}=0$ by $JT\in(P)$ *i.e.* Lemma 4.2.)

Let for convenience $\mathcal{N}_{\pm\mathrm{i}}=\operatorname{Ker}_{\pm\mathrm{i}}\dot{T}^*$, where $\dot{T}=JT$ in $(\mathcal{H},[\cdot,J\cdot])$.

Because $(J\pm\mu I)^{-1}\subseteq(1-\mu^2)^{-1}(J\mp\mu I)$, therefore

$$\operatorname{Ker}_\lambda\Sigma=(J-\mu I)(\mathcal{N}_{\mathrm{i}}\cap L_\eta(\mathcal{N}_{-\mathrm{i}}))$$

where the operator

$$L_\eta = \eta P^- + \frac{1}{\eta} P^+ , \quad \eta = \frac{1-\mu}{1+\mu} .$$

From here $\mathrm{Ker}_\lambda\,\Sigma$ is trivial iff so is $\mathcal{N}_{\mathrm{i}} \cap L_\eta(\mathcal{N}_{-\mathrm{i}})$, *i.e.* the system

$$(5.1) \qquad \begin{aligned} &(x^- + x^+, -\mathrm{i}x^- - \mathrm{i}x^+) \in \dot{T}^* , \\ &(\eta x^- + \eta^{-1}x^+, \mathrm{i}\eta x^- + \mathrm{i}\eta^{-1}x^+) \in \dot{T}^* \end{aligned}$$

has only trivial solutions for $x^\pm \in \mathcal{H}^\pm$.

From (5.1) follows

$$((\eta^2-1)x^-, \mathrm{i}(\eta^2+1)x^- + 2\mathrm{i}x^+) \in \dot{T}^* \,|\, \mathcal{H}^-$$

and then

$$(x^-, \zeta x^-) \in P^-\dot{T}^* \,|\, \mathcal{H}^- , \quad \zeta = \mathrm{i}\frac{\eta^2+1}{\eta^2-1} (= \frac{1-\lambda^2}{2\lambda}) .$$

Since $\mathrm{Im}\,\zeta = -\dfrac{1+|\lambda|^2}{2|\lambda|^2}\,\mathrm{Im}\,\lambda$, by Lemma 2.4 a) $x^- = 0$ and then by (5.1) and Lemma 4.2 $x^+ = 0$ for all $\lambda \in \mathbb{C}_* \smallsetminus \{-\mathrm{i}, \mathrm{i}\}$; *i.e.* $\sigma_p(\Sigma) \subseteq (\mathbb{R} \smallsetminus \{0\}) \sqcup \{\pm\mathrm{i}\}$.

The eigenspace $\mathrm{Ker}_{\mathrm{i}}\,\Sigma$ consists of $x_{\mathrm{i}} + x_{-\mathrm{i}}$, with $x_{\pm\mathrm{i}} \in \mathcal{N}_{\pm\mathrm{i}}$, $(J-I)x_{\mathrm{i}} = (J+I)x_{-\mathrm{i}}$, *i.e.* $x_{\mathrm{i}} \in \mathcal{H}^+$ and $x_{-\mathrm{i}} \in \mathcal{H}^-$; but $\mathcal{H}^- \cap \mathcal{N}_{-\mathrm{i}} = \{0\}$. Similarly one establishes $\mathrm{Ker}_{-\mathrm{i}}\,\Sigma$.

From (5.1) also follows

$$(\frac{\eta^2-1}{2\mathrm{i}}x^-, x^+) \in P^+\dot{T}^* \,|\, \mathcal{H}^- , \quad (x^+, \frac{2\mathrm{i}\eta^2}{1-\eta^2}x^-) \in P^-\dot{T}^* \,|\, \mathcal{H}^+$$

*i.e.*

$$(x^-, (\frac{2\eta}{1-\eta^2})^2 x^-) \in (P^-\dot{T}^* \,|\, \mathcal{H}^+)(P^+\dot{T}^* \,|\, \mathcal{H}^-) .$$

If $T \in (LP')$, $P^-\dot{T}^* \,|\, \mathcal{H}^+ \subseteq (P^+\dot{T}^* \,|\, \mathcal{H}^-)^*$ by Lemma 2.4 b). Since the eigenvalues of $(P^+\dot{T}^* \,|\, \mathcal{H}^-)^*(P^+\dot{T}^* \,|\, \mathcal{H}^-)$ are real non-negative, while

$$(\frac{2\eta}{1-\eta^2})^2 = -\frac{1}{4}(\lambda + \frac{1}{\lambda})^2 < 0$$

($\lambda \neq 0$ is real), so $x^- = 0$ and $x^+ = 0$; *i.e.* $\sigma_p(\Sigma) \subseteq \{\pm\mathrm{i}\}$. □

Therefore $\mathbb{C}_* \cap \sigma_p(N) \subseteq \{\pm\mathrm{i}\} \cap \sigma_p(JT_0)$, and if $T \in (LP')$ then moreover $\sigma_p(N) \subseteq \{\pm\mathrm{i}\} \cap \sigma_p(JT_0)$. This implies by Corollary 3.3 the following characterization of $\mathbb{C}_* \cap \sigma_p(T_0) \smallsetminus \sigma_p(JT_0)$: $\lambda \in \mathbb{C}_* \smallsetminus \sigma_p(JT_0)$ belongs to $\sigma_p(T_0)$ iff either $\lambda \in \sigma_p(T)$ or $\mathcal{R}_{T-\lambda I} \cap \mathcal{R}_{N-\lambda}$ is nontrivial; and similarly in case $T \in (LP')$ for $\sigma_p(T_0) \smallsetminus \sigma_p(JT_0)$.

*Example* 5.3. As a handy illustration let in $\mathcal{H} = \mathbb{C}^2$, $J(x,y) = (-x,y)$, $[(x_1,x_2)\,,\,(y_1,y_2)] = -\bar{x}_1 y_1 + \bar{x}_2 y_2$, the symmetric operator $T(x,0) = (0,-\beta x)$, $\mathcal{D}_T = \mathcal{H}^- = \mathbb{C}\times\{0\}$, $\beta > 0$, and let the self-adjoint extension $T_0(x,y) = (\beta y, -\beta x)$, $\mathcal{D}_{T_0} = \mathcal{H}$. Then $N(0,y) = (\beta y, 0)$, $\mathcal{D}_N = \mathcal{H}^+ = \{0\}\times\mathbb{C}$, $\sigma_p(N) = \emptyset = \sigma_p(T)$ (*cf.* $\sigma_p(\Sigma) = \{\pm 1\}$), so $\sigma_p(T_0) = \{\pm \mathrm{i}\beta\} = O^c(T,N)$. Indeed, $\mathcal{R}_{T-\lambda I} \cap \mathcal{R}_{N-\lambda I} = \{(\lambda x, \beta x)\,|\,x\in\mathbb{C}\}$ if $\lambda = \pm\mathrm{i}\beta$, and $= \{0\}$ if otherwise. Clearly $\sigma_p(T_0) \subset C_T = \{\lambda \mid |\mathrm{Im}\,\lambda| \leq \beta\}$. Moreover, because $\mathcal{L}_T = \{0\}$, $T \in (LP)$ but $T \notin (LP')$, and this reflects a general observation: $T \neq T^c$ of class $(LP')$ exists only in an infinite-dimensional $\mathcal{H}$; [38, Corollary 3.4], [42, Corollaries 4.3, 7.6] are relevant to this. (For relations whose domain and range are orthogonal we refer to [27].)

5.4 Can a $J$-symmetric relation be made into that of class $(LP)$? In order that a closed $J$-symmetric relation, $A$ say, should have $T \in \mathrm{Sym}(A)$ of class $(LP)$, clearly it is necessary that the $J$-adjoint $A^c \in (L)$; *cf.* $A = A^c \in (L)$ is of class $(LP')$.

5.5 Let $T$ be closed $J$-symmetric. The relation $T_c = T^c\,|\,\mathcal{D}_T$ in (2.1) is $J$-symmetric, $J$-self-adjoint iff $\mathcal{D}_T = \bar{\mathcal{D}}_T \cap \mathcal{D}_{T^c}$, but not closed in general. We refer to *e.g.* [23, 25] for a Hilbert space analogue.

If $T$ has equal defect numbers, that is, if there is an ordinary boundary triple (obt) $(\mathcal{G},\Gamma)$ for $T^c$, then [25] $T_c$ is 1-1 with the so-called (*e.g.* [14, 16]) forbidden relation

$$\mathcal{F}_\Gamma = \Gamma(\{0\}\times\mathcal{H})\,.$$

We remove the restriction imposed on defect numbers, and as a result we give a criterion in order that a symmetric relation should have a symmetric extension of class $(LP')$.

5.6 Let $S$ be a closed $J$-symmetric relation, $S^c$ its $J$-adjoint. A unitary boundary pair (ubp) $(\mathcal{G},\Gamma)$ for $S^c$—which exists always [10, 15]—is as in [10, Definition 3.1]; here $\mathcal{G}$ is a Hilbert (parameter) space and $\Gamma$ a boundary relation with $\bar{\mathcal{D}}_\Gamma = S^c$.

**Proposition 5.4.** ([31, Proposition 4.5]) *Let $\Gamma$ be a boundary relation for $S^c$ with the closed domain, $\mathcal{D}_\Gamma = S^c$. Then $T \in \mathrm{Ext}(S)$ is in 1-1 correspondence with $\theta \in \mathrm{Ext}(\mathrm{Ind}\,\Gamma)$, denoted $T \leftrightarrow \theta$, given by*

$$T = \Gamma^{-1}(\theta)\,, \quad \theta = \Gamma(T)\,.$$

*If $T \leftrightarrow \theta$ then $T$ is closed iff so is $\theta$, in which case $T^c \leftrightarrow \theta^*$.*

*Remark* 5.5. Not just a ubp exists always, but also a ubp with the closed domain: Use the actual construction in [15, Proposition 3.7] and then apply *e.g.* [10, Proposition 3.3], see also [31, Example 6.3].

**Corollary 5.6.** *Let* $S_c = S \hat{+} (\{0\} \times \operatorname{Ind} S^c)$. *Then*

$$\operatorname{Sym}(S_c) = \{T \in \operatorname{Sym}(S) \mid \operatorname{Ind} T = \operatorname{Ind} S^c\} \leftrightarrow \operatorname{Sym}(\mathcal{F}_\Gamma)\,.$$

(Here $\leftrightarrow$ for sets is self-explanatory.)

*Proof.* Clearly $S_c \leftrightarrow \mathcal{F}_\Gamma$. Let $T \in \operatorname{Sym}(S)$, $\operatorname{Ind} T = \operatorname{Ind} S^c$, $T \leftrightarrow \theta \in \operatorname{Ext}(\operatorname{Ind}\Gamma)$. Then $\theta \in \operatorname{Sym}(\mathcal{F}_\Gamma)$ since $\theta \in \operatorname{Sym}(\operatorname{Ind}\Gamma)$ and

$$\theta = \Gamma(T) \supseteq \Gamma(\{0\} \times \operatorname{Ind} T) = \Gamma(\{0\} \times \operatorname{Ind} S^c) = \mathcal{F}_\Gamma\,. \qquad \square$$

(If $\operatorname{Sym}(S) \ni T \leftrightarrow \theta \in \operatorname{Sym}(\operatorname{Ind}\Gamma)$, in order that $\operatorname{Ind} T = \operatorname{Ind} T^c$ it is necessary that $\theta \cap \mathcal{F}_\Gamma = \theta^* \cap \mathcal{F}_\Gamma$, but not sufficient.)

A ubp $(\mathcal{G}, \Gamma)$ for $S^c$ is an $S$-generalized boundary pair (*cf.* [18, Definition 5.11]) if $\Gamma^{-1}(\{0\} \times \mathcal{G})$ is $J$-self-adjoint. This holds if *e.g.* $(\mathcal{G}, \Gamma)$ is an obt for $S^c$, *i.e.* a ubp with surjective $\Gamma$; the latter needs $S$ to have equal defect numbers, see *e.g.* [10] for further discussion.

**Proposition 5.7.** *Let $A$ be a closed symmetric relation in a Krein space $\mathcal{H} = \mathcal{H}^-[\oplus]\mathcal{H}^+$; $A^c = L \hat{+} (\{0\} \times \operatorname{Ind} A^c)$ for some closed operator $L$ (e.g. the operator part of $A^c$). If $\mathcal{H}^- \subseteq \mathcal{D}_L$ and the operator $P^- L \mid \mathcal{H}^-$ is self-adjoint in $(\mathcal{H}^-, -[\cdot\,,\cdot])$, then there is $T \in \operatorname{Sym}(A)$ of class $(LP')$.*

*Proof.* Let $(\mathcal{G}, \dot{\Gamma})$ be a ubp for $A^c = \mathcal{D}_{\dot{\Gamma}}$; such a ubp exists as pointed out in Remark 5.5. Consider a closed relation $H \subseteq A^c \cap H^c$ in $\mathcal{H}$. Put

$$S = A \hat{+} H \leftrightarrow \theta_H = \dot{\Gamma}(H)\,.$$

Then in particular $S \in \operatorname{Sym}(A)$. The pair $(\mathcal{G}, \Gamma)$ where

$$\Gamma = (\dot{\Gamma} \mid H^c) \hat{+} (H \times \{0\})$$

(with domain $S^c$, kernel $S$, and multivalued part $\theta_H$) is a ubp for $S^c$; one might prefer $\Gamma = \dot{\Gamma} V_H$, where $V_H = (I \mid H^c) \hat{+} (\{0\} \times H)$ is self-adjoint, unitary in a $J_\Gamma$-space. (If $\mathcal{D}_{\theta_H} = \overline{\mathcal{D}}_{\theta_H}$, $(\mathcal{G}, \Gamma)$ is an $S$-generalized boundary pair.)

If $\mathcal{H}^- \subseteq \mathcal{D}_H$ then $\operatorname{Sym}(S_c) \ni T \leftrightarrow \theta \in \operatorname{Sym}(\mathcal{F}_\Gamma)$ (Corollary 5.6) is of class $(LP')$. We construct such a $H$.

Define a bounded operator $H$ ($\mathcal{D}_H = \mathcal{H}^-$) in $\mathcal{H}$ by

$$H = L \,|\, \mathcal{H}^- = L^{--} + L^{+-}\,,$$

$$L^{--} = P^- L \,|\, \mathcal{H}^-\,, \quad L^{+-} = P^+ L \,|\, \mathcal{H}^-\,.$$

Regarding $L^{--}$ and $L^{+-}$ as Hilbert space operators $\mathcal{H}^- \to \mathcal{H}^-$ and $\mathcal{H}^- \to \mathcal{H}^+$, both with domain $\mathcal{H}^-$, the $J$-adjoint

$$H^c = L^{--} \hat{+} -(L^{+-})^* \hat{+} (\{0\} \times \mathcal{H}^+)$$

hence $H \subseteq H^c \,|\, \mathcal{H}^- \subseteq H^c$ and in fact $H \subseteq A^c \cap H^c$. $\square$

*Remark* 5.8. Consider the operator part $A_s$ of $A$ as an operator in the Krein space $(\bar{\mathscr{D}}_{A^c}, [\cdot\,,\,\cdot])$, $\mathscr{H}^- \subseteq \mathscr{D}_{A^c}$; $A_s^c = (A^c)_s \mathbin{\widehat{\oplus}} (\{0\}\times(\mathscr{D}_A + \operatorname{Ind} A)^{[\perp]})$ is its $J$-adjoint, where $(A^c)_s$ is the operator part of $A^c$. Let $A_s^{--} = P^- A_s \,|\, \mathscr{H}^-$ in $(\mathscr{H}^-, -[\cdot\,,\,\cdot])$; the adjoint $(A_s^{--})^* = \overline{P^- A_s^c \,|\, \mathscr{H}^-}$. From here, if $\operatorname{Ind} A = \operatorname{Ind} A^c$, $P^- L \,|\, \mathscr{H}^- = A_s^{--}$ is self-adjoint in the Hilbert space $\mathscr{H}^-$, *cf.* [33]. The condition is not necessary though; as an illustration, take $H = k^{-1} I^-$, $k > \pi^{-2}$, for $A$ in [47, Section 3] to deduce that a non-densely defined and non-maximal symmetric $T = S_c \in \operatorname{Sym}(A)$ is of class $(LP')$.

## 6. On points of regular type

6.1 Let $T$ be a closed relation in a $J$-space. Once $\lambda_0$ such that $\mathscr{R}_{T-\lambda_0 I}$ (or equivalently $\mathscr{R}_{T^c-\bar{\lambda}_0 I}$) is closed has been established, $\mathscr{R}_{T-\lambda I}$ is closed for $\lambda$ in a neighborhood of $\lambda_0$ [36, Lemma 2.7]. Here we look for $\lambda_0 \in r(JT)$.

6.2 Let $T$ be a relation in a $J$-space with canonical decomposition $\mathscr{H} = \mathscr{H}^-[\oplus]\mathscr{H}^+$. One verifies (*cf.* Section 2.5; $I^+ = I \,|\, \mathscr{H}^+$)

$$\begin{aligned} P^+(JT + \lambda I)^{-1} \,|\, \mathscr{H}^+ &= (\mathscr{T}^+(\lambda) + \lambda I^+)^{-1}\,, \\ P^+(T - \lambda I)^{-1} \,|\, \mathscr{H}^+ &= (\mathscr{T}^+(\lambda) - \lambda I^+)^{-1} \end{aligned} \tag{6.1}$$

for all $\lambda \in \mathbb{C}$, where the relation $\mathscr{T}^+(\lambda)$ in $\mathscr{H}^+$ is defined by

$$\mathscr{T}^+(\lambda) = \{(P^+x, P^+y) \,|\, (x,y) \in T\,;\; y - \lambda x \in \mathscr{H}^+\}\,.$$

**Lemma 6.1.** $\mathscr{R}_{\mathscr{T}^+(\lambda)-\lambda I^+}$ $(\lambda \in \mathbb{C})$ *equals each of the following lineals:*

a) $\mathscr{H}^+ \cap \mathscr{R}_{T-\lambda I}$.

b) $P^+(I - 2\lambda(JT + \lambda I)^{-1})(\mathscr{H}^+)$*;* $(I - 2\lambda(JT + \lambda I)^{-1})(\mathscr{H}^+)$ *is a subspace for* $\lambda \neq 0$*, if* $T$ *is closed and* $\mathscr{H}^- + \mathscr{R}_{JT^c-\bar{\lambda}I}$ *is a subspace.*

c) $P^+(JT - \lambda I)(JT + \lambda I)^{-1}(\mathscr{H}^+)$ *if* $\operatorname{Ind} T \subseteq \mathscr{H}^+$.

*Proof.* a) Clear from the definition of $\mathscr{T}^+$.

b) Since the first statement is clear, we only verify the second one, and so let $T = \bar{T}$, $V_\lambda = I + 2\lambda(JT - \lambda I)^{-1}$, $\lambda \neq 0$. Then $E_{\bar{\lambda}} = 2\bar{\lambda}(JT^c - \bar{\lambda} I)^{-1}$ is a closed relation in $\mathscr{H}$, with the (Hilbert) adjoint $E_{\bar{\lambda}}^* = 2\lambda(JT - \lambda I)^{-1}$. By Lemma 2.1 a) $\mathscr{H} \ominus \mathscr{H}_{\bar{\lambda}} = V_\lambda^{-1}(\mathscr{H}^+)$ for the lineal $\mathscr{H}_{\bar{\lambda}} = (I + E_{\bar{\lambda}})(\mathscr{H}^-)$ if $\mathscr{H}^- + \mathscr{R}_{JT^c-\bar{\lambda}I}$ is a subspace; $V_\lambda^{-1} = V_{-\lambda}$ ([12, Lemma 1.1.9]).

c) By definition $\mathscr{R}_{\mathscr{T}^+(\lambda)-\lambda I^+}$ is the set of $P^+y$ such that $(\exists x)$ $(x,y) \in JT - \lambda I$, $y + 2\lambda x \in \mathscr{H}^+$. On the other hand, the relation $(JT - \lambda I)(JT + \lambda I)^{-1}$ consists of $(v, y) \in \mathscr{H}_\Gamma$ such that $(\exists x)$ $(x,y) \in JT - \lambda I$, $(x, v) \in JT + \lambda I$; equivalently, such that $(\exists x)$ $(x,y) \in JT - \lambda I$, $y - v + 2\lambda x \in \operatorname{Ind} JT$. Therefore $P^+(JT - \lambda I)(JT + \lambda I)^{-1}(\mathscr{H}^+)$ is the set of $P^+y$ such that $(\exists x)$ $(x,y) \in JT - \lambda I$, $y + 2\lambda x \in \mathscr{H}^+ + \operatorname{Ind} JT$. □

**Lemma 6.2.** *Let* $T$ *be a closed relation in a Krein space with canonical decomposition* $\mathscr{H} = \mathscr{H}^-[\oplus]\mathscr{H}^+$ *and canonical symmetry* $J$*. Fix a nonzero*

*$\lambda \in \mathbb{C}$ such that $\mathscr{H}^- + \mathscr{R}_{JT^c - \bar\lambda I}$ is a subspace. Then $\mathscr{R}_{T-\lambda I}$ is a subspace iff so is $\mathscr{H}^- \cap \mathscr{R}_{T^c - \bar\lambda I}$; if also $\mathscr{H}^+ + \mathscr{R}_{JT+\lambda I}$ is a subspace, then $\mathscr{R}_{T-\lambda I}$ is a subspace iff so is $\mathscr{H}^+ \cap \mathscr{R}_{T-\lambda I}$.*

*Proof.* By direct inspection

$$\mathscr{R}_{T-\lambda I} = J(\mathscr{R}_{I_\lambda}), \quad I_\lambda = I + 2\lambda P^-(JT - \lambda I)^{-1}.$$

for all $\lambda \in \mathbb{C}$. The adjoint $I_{\bar\lambda}^*$ in $(\mathscr{H}, [\cdot, J \cdot])$ consists of $(y, x + y) \in \mathscr{H}_\Gamma$ with $(x, \bar\lambda(I - J)y) \in JT^c - \bar\lambda I$; hence

$$\mathscr{R}_{I_{\bar\lambda}^*} = (\mathscr{H}^- \cap \mathscr{R}_{T^c - \bar\lambda I})[\oplus]\mathscr{H}^+ \quad (\lambda \neq 0).$$

Now, the proof of the first iff statement is accomplished after we show that $I_\lambda$ $(\lambda \neq 0)$ is a closed relation in $\mathscr{H}$ if $\mathscr{H}^- + \mathscr{R}_{JT^c - \bar\lambda I}$ is a subspace.

Let $E_{\bar\lambda}$ be a closed relation as in the proof of Lemma 6.1 b). Then $I_\lambda = I + P^- E_{\bar\lambda}^*$ is a closed relation in $\mathscr{H}$ if such is $P^- E_{\bar\lambda}^*$. But the latter relation in $\mathscr{H}$ is closed if $\mathscr{H}^- + \mathscr{D}_{E_{\bar\lambda}}$ $(= \mathscr{H}^- + \mathscr{R}_{JT^c - \bar\lambda I})$ is a subspace.

By the above, in order to prove the second iff statement it suffices to show for a nonzero $\lambda$ that $\mathscr{H}^- \cap \mathscr{R}_{T^c - \bar\lambda I}$ is a subspace if such are $\mathscr{H}^- + \mathscr{R}_{JT^c - \bar\lambda I}$, $\mathscr{H}^+ + \mathscr{R}_{JT+\lambda I}$, and $\mathscr{H}^+ \cap \mathscr{R}_{T-\lambda I}$.

If $\mathscr{H}^- + \mathscr{R}_{JT^c - \bar\lambda}$ and $\mathscr{H}^+ + \mathscr{R}_{I+E_{\bar\lambda}^*}$ $(= \mathscr{H}^+ + \mathscr{R}_{JT+\lambda I})$ are subspaces, by Lemma 2.1 b) the lineal $\mathscr{H}_{\bar\lambda}$ in the proof of Lemma 6.1 b) is a subspace. For a subspace $\mathscr{H}_{\bar\lambda}$, $P^-(\mathscr{H}_{\bar\lambda})$ is a subspace iff such is $P^+(\mathscr{H}_{\bar\lambda}^\perp)$ (Lemma 2.1 c)). Now recall that $P^+(\mathscr{H}_{\bar\lambda}^\perp) = \mathscr{H}^+ \cap \mathscr{R}_{T-\lambda I}$, and similarly one verifies the equality $P^-(\mathscr{H}_{\bar\lambda}) = \mathscr{H}^- \cap \mathscr{R}_{T^c - \bar\lambda I}$. □

If $\mathscr{H}^+ \cap \mathscr{R}_{JT+\lambda I}$ is trivial for some $\lambda \in \mathbb{C} \setminus (\{0\} \cup \sigma_p(JT) \cup -\sigma_p(JT))$ then so is $\mathscr{H}^+ \cap \mathscr{R}_{T-\lambda I}$:

*Example* 6.3. a) Let $T$ be a closed relation in $\mathscr{H}$ such that $\mathscr{H}^- + \mathscr{R}_{JT^c - \bar\lambda I}$ and $\mathscr{H}^+ + \mathscr{R}_{JT+\lambda I}$ are subspaces and moreover $\mathscr{H}^+ \cap \mathscr{R}_{JT+\lambda I} = \{0\}$, some $\lambda \in \mathbb{C} \setminus (\{0\} \cup \sigma_p(JT) \cup -\sigma_p(JT))$. Then $\mathscr{R}_{T-\lambda I}$ is a subspace and $\mathscr{H}^+ \cap \mathscr{R}_{T-\lambda I} = \{0\}$.

*Proof.* Straightforwardly

$$P^+(\mathrm{Ker}_{-\lambda}(JT \hat{+} \lambda I^+)) = \mathscr{H}^+ \cap \mathscr{R}_{T-\lambda I}$$

for $\lambda \in \mathbb{C} \setminus \{0\}$. By Proposition 3.2 $\mathrm{Ker}_{-\lambda}(JT \hat{+} \lambda I^+) = \mathrm{Ker}_{-\lambda} JT$ iff $\mathscr{H}^+ \cap \mathscr{R}_{JT+\lambda I} = \mathscr{H}^+ \cap \mathrm{Ker}_\lambda JT$. Then apply Lemma 6.2. □

b) As a simple illustration let $T = J$ on a maximal negative subspace.

*Remark* 6.4. Let $T$ be a relation in a $J$-space. Then

$$\mathscr{H}^+ + \mathscr{R}_{JT+\lambda I} = \mathscr{H}^+ + \mathscr{R}_{T-\lambda I} \quad \text{and}$$
$$\mathscr{H}^- + \mathscr{R}_{JT-\lambda I} = \mathscr{H}^- + \mathscr{R}_{T-\lambda I} \quad (\lambda \in \mathbb{C}).$$

*Proof.* Let $y^- \in P^-(\mathcal{R}_{T-\lambda I})$, then $y^- = P^- y$ and $(\exists x)$ $(x, y) \in T - \lambda I$, so $(x, Jy + 2\lambda P^+ x) \in JT + \lambda I$. With $v = Jy + 2\lambda P^+ x$ one concludes that $y^- = -P^- v$, $(x, v) \in \mathcal{R}_{JT+\lambda I}$, and since the argument is reversible this proves the first equality; the proof of the second one is analogous. □

6.3 Let $T$ be a closed relation in a $J$-space with canonical decomposition $\mathcal{H} = \mathcal{H}^-[\oplus]\mathcal{H}^+$, and let $s(JT)$ denote the union of $\{0\} \cap r(JT)$ with the set of $\lambda \in r(JT) \setminus \{0\}$ such that $\mathcal{H}^+ + \mathcal{R}_{JT+\lambda I}$ and $P^+(I - 2\lambda(JT + \lambda I)^{-1})(\mathcal{H}^+)$ are subspaces.

By Lemmas 6.1, 6.2 $s(JT) \subseteq \{\lambda \,|\, \mathcal{R}_{T-\lambda I} = \bar{\mathcal{R}}_{T-\lambda I}\}$.

Some corollaries and comments follow.

6.4 $s(JT) \cap \hat{\sigma}(T) = s(JT) \cap \sigma_p(T)$.

6.5 $\lambda \in r(JT) \setminus \{0\}$ is a $\Phi_+$-point of a closed relation $T$ in a $\Pi_\kappa$-space [21, 35], *i.e.* $\dim \mathrm{Ker}_\lambda T = \dim (JT - \lambda I)^{-1}(\mathrm{Ker}_\lambda \mathcal{T}^-(\lambda)) \leq \kappa$, $\mathcal{R}_{T-\lambda I} = \bar{\mathcal{R}}_{T-\lambda I}$ (because $s(JT) = r(JT)$ by Lemma 6.1 b) and recall the proof of Lemma 2.7 d)). If what is more $T \in (L)$, a number of results adapted from [21, Section 8] can be formulated, which, however, we skip here.

6.6 Proof of Theorem 1.3 c): Because $\mathbb{C}^- \subseteq r(JT_0)$, $\Delta_0^- = \mathbb{C}^- \cap s(JT_0)$ by Lemma 6.1 c). Thus the assertion follows from a) and Sections 6.4, 6.5.

## 7. A FEW REMARKS ON A SELF-ADJOINT RELATION

7.1 For a $J$-self-adjoint relation $T$ the inclusion $\mathbb{C}_* \cap s(JT) \subseteq s_*(T) = \{\lambda \in \mathbb{C}_* \,|\, \mathcal{R}_{T-\lambda I} = \bar{\mathcal{R}}_{T-\lambda I}\}$ becomes the equality, the set $s_*(T)$ is symmetric with respect to the real axis (in symbols $s_*^*(T) = \{\lambda \,|\, \bar{\lambda} \in s_*(T)\} = s_*(T)$), and

$$\begin{aligned} \mathbb{C}_* \cap r(T) &= s_*(T) \cap r(T) = s_*(T) \setminus \sigma_p(T)\,, \\ \mathbb{C}_* \cap \rho(T) &= s_*(T) \cap \rho(T) = s_*(T) \cap \delta(T) \end{aligned} \tag{7.1}$$

where $\delta(T) = (\mathbb{C}_* \setminus \sigma_p(T)) \cap (\mathbb{C}_* \setminus \sigma_p(T))^*$ as in [30, Definition 4.1].

Contrary to $s_*(T) = \mathbb{C}_*$ for a $\pi$-self-adjoint relation there is a $J$-self-adjoint relation with $s_*(T) = \mathbb{C}_*$, [19, Proposition 4.5], as well as $J$-self-adjoint operator with $s_*(T) = \emptyset$, [6, Example 2.3.31].

7.2 Let $T$ be a $J$-self-adjoint relation. Although $\sigma(T)$ is symmetric with respect to the real axis, that does not necessarily apply to eigenvalues. Given $s = s^* \subseteq s_*(T)$, $s \cap \sigma(T) \neq \emptyset$, the last formula in (7.1) shows that the part in $s$ of $\sigma_p(T)$ is symmetric with respect to the real axis iff $T$ has no residual spectrum, *i.e.* $\lambda \in \sigma(T) \setminus \sigma_p(T)$ with $\bar{\mathcal{R}}_{T-\lambda I} \neq \mathcal{H}$ ([6, Section 2.6]), in $s$.

Such $s$ for a definitizable relation $T$ consists of only non-real $\Phi_0$-points: $\lambda \in s_*(T)$ such that $\dim \mathrm{Ker}_\lambda T = \dim \mathrm{Ker}_{\bar{\lambda}} T < \infty$ [19, Proposition 5.2]; for an operator in a $\Pi_\kappa$-space the result goes back to [33, 39].

Such $s$ for $T \in (L)$ consists of only $\Phi_0$-points too (Proposition 7.1).

7.3 As already remarked in Section 5.4 a $J$-self-adjoint relation $T \in (L)$ is automatically of class $(LP')$. As a special case, a $\pi$-self-adjoint operator $T$ is of class $(LP')$, with spectrum contained in $\{\lambda \mid |\mathrm{Im}\,\lambda| \leq \|P^+TP^-\|\}$ $(\Pi^- \subseteq \mathcal{D}_T)$ [33, Theorem 2.1]. By the same perturbation argument [21, Sections 5.1, 5.3] a similar estimate holds for a $J$-self-adjoint relation $T \in (LP')$.

**Proposition 7.1.** *Let $T \in (L)$ be a self-adjoint relation in a Krein space $(\mathcal{H}, [\cdot, \cdot])$ with canonical decomposition $\mathcal{H} = \mathcal{H}^-[\oplus]\mathcal{H}^+$; $\mathcal{H}^- \subseteq \mathcal{D}_T$. Then*

$$\sigma(T) \subseteq \{\lambda \mid |\lambda - \lambda_0| \leq \|P^+T_sP^-\| ;\\ \lambda_0 \in \sigma(P^-T \,|\, \mathcal{H}^-) \cup \sigma(P^+T_s \,|\, \mathcal{H}^+)\}\,.$$

*Particularly*

a) *$\mathbb{C}_* \cap \sigma(T)$ is the set of non-real normal eigenvalues of $T_s$.*

b) *$s_*(T)$ is a nonempty subset of normal points of $T$.*

As usual ([6, Definition 2.1.14]) a normal point is either a normal eigenvalue or a regular point.

*Proof.* With the notation as in Section 5.2 and by Lemma 2.4 a self-adjoint operator $T_s = H + B$ in a Krein space $\mathcal{H}_c$ is the sum of a diagonal (with respect to $\mathcal{H}_c = \mathcal{H}^-[\oplus]\mathcal{H}_c^+$) self-adjoint operator

$$H = \begin{pmatrix} T_s^{--} & 0 \\ 0 & T_s^{++} \end{pmatrix}$$

both in a Krein space $(\mathcal{H}_c, [\cdot, \cdot])$ and in a Hilbert space $(\mathcal{H}_c, [\cdot, J\cdot])$ $(J = P^+ - P^-)$, and an off-diagonal bounded operator

$$B = \begin{pmatrix} 0 & -(T_s^{+-})^* \\ T_s^{+-} & 0 \end{pmatrix}$$

in $\mathcal{H}_c$. Because ([21, p. 76]) each $\lambda \in \sigma(T_s)$ is contained in a disc $|\lambda - \lambda_0| \leq \|B\|$ for some $\lambda_0 \in \sigma(H)$, and because $\|B\| \leq \|T_s^{+-}\|$ and $\sigma(T) = \sigma(T_s)$ (Section 5.1), the first assertion follows.

a) This is an application of [6, Corollary 2.2.12, Proposition 2.3.7].

b) If $s_*(T) = \emptyset$ then by (7.1) $\mathbb{C}_* \subseteq \sigma(T)$—a contradiction to what has been obtained. By a), that $\lambda \in s_*(T)$ is a normal point is just another way of saying that $\lambda$ is not a point from the residual spectrum of $T$. □

7.4 Consider a relation $T = T^c$ in a Krein space $(\mathcal{H}, [\cdot, \cdot])$ with canonical decomposition $\mathcal{H} = \mathcal{H}^-[\oplus]\mathcal{H}^+$ and canonical symmetry $J$. Here one finds an interpretation of the two equalities

$$\begin{aligned} \mathbb{C}_* \cap \rho(T) &= \{\lambda \in \mathbb{C}_* \,|\, 0 \in \rho(\mathcal{T}^-(\lambda) - \lambda I^-)\} \\ &= \{\lambda \in \mathbb{C}_* \,|\, 0 \in \rho(\mathcal{T}^+(\lambda) - \lambda I^+)\} \end{aligned} \tag{7.2}$$

by means of the resolvent formula of Krein–Naimark type for a Nevanlinna function $\mathcal{T}^-$ (Lemma 2.7 a), b)) in $(\mathcal{H}^-, -[\cdot,\cdot])$ via a Nevanlinna family $\mathcal{T}^+$ ($\operatorname{Ind}\mathcal{T}^+(\lambda) = \operatorname{Ind} T$) in $(\mathcal{H}^+, [\cdot,\cdot])$ in case $\mathcal{H}^- \subseteq \mathcal{D}_T$ (*i.e.* $T \in (L)$); that $\mathcal{T}^+$ is Nevanlinna one shows by exactly the same argument as for $\mathcal{T}^-$.

The first equality in (7.2), in an arbitrary canonical decomposition, follows from Lemma 2.7 d), Remark 6.4, and $\mathcal{R}_{\mathcal{T}^-(\lambda)-\lambda I^-} = \mathcal{H}^- \cap \mathcal{R}_{T-\lambda I}$ (*cf.* Lemma 6.1 a)); the second one is shown analogously.

Nevanlinna families $\mathcal{T}^\pm$ are realized via Weyl families of boundary relations, [15], [12, Section 2.7], whose so-called main transform is $T$. Particularly $-\mathcal{T}^+$ is the Weyl family associated with a boundary relation of a symmetric relation $T \cap (\mathcal{H}^-)^2$ in $(\mathcal{H}^-, [\cdot,\cdot])$, and one recovers a familiar equality in (6.1), second formula, if $\lambda \in \rho(T)$.

**Proposition 7.2.** *Let* $T = T^c \in (L)$*,* $\mathcal{H}^- \subseteq \mathcal{D}_T$*,* $T^{--} = P^- T \,|\, \mathcal{H}^-$*, and* $\mathcal{L}_\lambda = (T-\lambda I)^{-1}(\mathcal{H}^+)$*,* $\lambda \in \mathbb{C}_*$*. Then*

a) $\mathcal{L}_\lambda$ *is of class* $\mathcal{A}^+$*, i.e.* $\mathcal{H}^- \cap \mathcal{L}_\lambda = \{0\}$*, and a closable operator*

$$\gamma_\lambda = P^-(P^+ \,|\, \mathcal{L}_\lambda)^{-1} \colon \mathcal{H}^+ \to \mathcal{H}^-\,,$$

$$\mathcal{D}_{\gamma_\lambda} = P^+(\mathcal{L}_\lambda) = \mathcal{H}^+ \cap \mathcal{D}_T = \mathcal{D}_{\mathcal{T}^+(\lambda)} \quad (\mathcal{D}_{\gamma_\lambda^c} = \mathcal{H}^-)$$

*is bounded;* $\gamma_\lambda$ *is therefore a contraction iff* $\mathcal{L}_\lambda$ *is non-negative.*

b) $(\mathcal{T}^-(\lambda) - \lambda I^-)^{-1} = (T^{--} - \lambda I^-)^{-1} + \gamma_\lambda(\mathcal{T}^+(\lambda) - \lambda I^+)^{-1}\gamma^c_{\bar\lambda}$.

*Proof.* a) That $\mathcal{L}_\lambda \in \mathcal{A}^+$ is seen from

$$\mathcal{L}_\lambda = \{x \mid (\exists y)\ (x,y) \in T\,;\ y - \lambda x \in \mathcal{H}^+\}$$

by applying Lemma 2.4 a).

The relation

$$\Gamma = \{((P^-x, P^-y), (P^+x, -P^+y)) \mid (x,y) \in T\}$$

from a $J_\Gamma$-space $((\mathcal{H}^-)^2, [\cdot,\cdot]_\Gamma)$ to a $J_\Gamma$-space $((\mathcal{H}^+)^2, [\cdot,\cdot]_\Gamma)$ is unitary, and using Lemma 2.4 a) $(\mathcal{H}^+, \Gamma)$ is an $S$-generalized boundary pair (Section 5.6) of a closed symmetric relation $\operatorname{Ker}\Gamma = T \cap (\mathcal{H}^-)^2$ in $(\mathcal{H}^-, [\cdot,\cdot])$; hence all the remaining arguments are essentially the same as in [18, Theorem 5.17], with $\gamma_\lambda$ (resp. $-\mathcal{T}^+$) the $\gamma$-field (resp. Weyl family) corresponding to $(\mathcal{H}^+, \Gamma)$.

b) This is a routine computation (*cf.* [18, Theorem 5.8]) by using that the adjoint $\gamma^c_{\bar\lambda} = -P^+(P^- \,|\, \mathcal{L}^{[\perp]}_{\bar\lambda})^{-1}$ where $\mathcal{L}^{[\perp]}_{\bar\lambda} = (T - \lambda I)(\mathcal{H}^-)$. □

*Remarks* 7.3. a) If $T$ is an operator, one recognizes $\mathcal{T}^+(\lambda) - \lambda I^+$ as the Schur complement $T(\lambda) = T^{++} - \lambda I^+ - T^{+-}(T^{--} - \lambda I^-)^{-1}T^{-+}$ (in the notation of Section 5.2) of $T - \lambda I$, so in particular by the resolvent formula $\mathcal{R}_{T(\lambda)} = \mathcal{H}^+$ iff $\mathcal{R}_{T(\lambda)} (= \mathcal{H}^+ \cap \mathcal{R}_{T-\lambda I}) \supseteq \mathcal{R}_{T^{+-}}$.

b) An interpretation of $-\mathscr{T}^+$ as Weyl family leads to a formula of well-known type: $\mathscr{T}^+(\lambda_0)-\mathscr{T}^+(\lambda) = (\lambda-\lambda_0)\gamma_{\bar\lambda}^c\gamma_{\lambda_0}$ $(\lambda_0, \lambda \in \mathbb{C}_*; \lambda_0 \neq \lambda)$.

## REFERENCES

[1] Y. Arlinskii and C. Tretter, *Everything is possible for the domain intersection* $\operatorname{dom} T \cap \operatorname{dom} T^*$, Adv. Math. **374** (2020), no. 18.

[2] Yu. Arlinskii, *Families of symmetric operators with trivial domains of their squares*, Compl. Anal. Oper. Theory **17** (2023), no. 120.

[3] ______, *Compressions of selfadjoint and maximal dissipative extensions of non-densely defined symmetric operators*, arXiv:2409.10234 (2024).

[4] ______, *Squares of symmetric operators*, Compl. Anal. Oper. Theory **18** (2024), no. 161.

[5] T. Azizov and I. Iokhvidov, *Linear operators in spaces with an indefinite metric and their applications*, Itogi nauki i mech. Ser. Mat. (in Russian) **17** (1979), 113–205.

[6] ______, *Linear Operators in Spaces with an Indefinite Metric*, John Wiley & Sons. Inc., 1989.

[7] T. Ya. Azizov, A. Dijksma, and I. V. Gridneva, *On the boundedness of Hamiltonian operators*, Proc. Amer. Math. Soc. **131** (2002), no. 2, 563–576.

[8] T. Ya. Azizov, A. Dijksma, and G. Wanjala, *Compressions of maximal dissipative and self-adjoint linear relations and of dilations*, Lin. Alg. Appl. **439** (2013), 771–792.

[9] T. Ya. Azizov and I. V. Gridneva, *On invariant subspaces of $J$-dissipative operators*, Ukrainian Math. Bull. **1** (2009), 1–13.

[10] J. Behrndt, V. A. Derkach, S. Hassi, and H. de Snoo, *A realization theorem for generalized Nevanlinna families*, Operators and Matrices **5** (2011), no. 4, 679–706.

[11] J. Behrndt, S. Hassi, and H. de Snoo, *Functional models for Nevanlinna families*, Opuscula Math. **28** (2008), no. 3, 233–245.

[12] ______, *Boundary Value Problems, Weyl Functions, and Differential Operators.*, Vol. 108, Birkhauser, 2020.

[13] J. Bognár, *Indefinite inner product spaces*, Springer-Verlag Berlin Heidelberg New York, 1974.

[14] V. Derkach, *On generalized resolvents of Hermitian relations in Krein spaces*, J. Math. Sci. **97** (1999), no. 5, 4420–4460.

[15] V. Derkach, S. Hassi, M. Malamud, and H. de Snoo, *Boundary relations and their Weyl families*, Trans. Amer. Math. Soc. **358** (2006), no. 12, 5351–5400.

[16] V. Derkach and M. Malamud, *The extension theory of Hermitian operators and the moment problem*, J. Math. Sci. **73** (1995), no. 2, 141–242.

[17] V. A. Derkach and M. M. Malamud, *Extension theory of symmetric operators and boundary value problems*, Vol. 104, Institute of Mathematics of NAS of Ukraine, Kiev, 2017 (in Russian).

[18] V. Derkach, S. Hassi, and M. M. Malamud, *Generalized boundary triples, I. Some classes of isometric and unitary boundary pairs and realization problems for subclasses of Nevanlinna functions*, Math. Nachr. **293** (2020), no. 7, 1278–1327.

[19] A. Dijksma and H. de Snoo, *Symmetric and self-adjoint relations in Krein spaces II*, Annales Academire Scientiarum Fennicae **12** (1987), 199–216.

[20] N. Dunford and J. T. Schwartz, *Linear Operators I: General Theory*, John Wiley and Sons, Inc., New York, 1988.

[21] I. Gokhberg and M. Krein, *Fundamental aspects of defect numbers, root numbers and indexes of linear operators*, Uspekhi Mat. Nauk (in Russian) **12** (1957), no. 2, 43–118.

[22] S. Hassi and H. de Snoo, *One-dimensional graph perturbations of selfadjoint relations*, Annales Academire Scientiarum Fennicae Mathematica **22** (1997), 123–164.

[23] S. Hassi and H. S. V. de Snoo, *Friedrichs and Krein type extensions in terms of representing maps*, Ann. Funct. Anal. **15** (2024), no. 78.

[24] S. Hassi, H. S. V. de Snoo, and F. H. Szafraniec, *Componentwise and Cartesian decompositions of linear relations*, Dissert. Math. **465** (2009), 1–59.

[25] ______, *Infinite-dimensional perturbations, maximally nondensely defined symmetric operators, and some matrix representations*, Indagationes Mathematicae **23** (2012), no. 4, 1087–1117.

[26] S. Hassi, J. Labrousse, and H. de Snoo, *Operational calculus for rows, columns, and blocks of linear relations*, Advances in Operator Theory **5** (2020), no. 3, 1193–1228.

[27] ______, *Selfadjoint extensions of relations whose domain and range are orthogonal*, Methods Func. Anal. Topology **26** (2020), no. 1, 39–62.

[28] S. Hassi, Z. Sebestyén, H. S. V. de Snoo, and F. H. Szafraniec, *A canonical decomposition for linear operators and linear relations*, Acta Math. Hungar. **115** (2007), no. 4, 281–307.

[29] I. Iokhvidov and M. Krein, *Spectral theory of operators in space with indefinite metric. I*, Tr. Mosk. Mat. Obs. (in Russian) **5** (1956), 367–432.

[30] R. Juršėnas, *On the similarity of boundary triples of symmetric operators in Krein spaces*, Compl. Anal. Oper. Theory **17** (2023), no. 72, 1–39.

[31] ______, *Weyl families of transformed boundary pairs*, Math. Nachr. **296** (2023), no. 8, 3411–3448.

[32] M. Krasnoselskii, *On self-adjoint extensions of Hermitian operators*, Ukrainian Math. J. (in Russian) **1** (1949), 21–38.

[33] M. Krein and H. Langer, *Defect subspaces and generalized resolvents of an Hermitian operator in the space* $\Pi_\kappa$, Funktsional. Anal. i Prilozhen. (in Russian) **5** (1971), no. 2, 59–71.

[34] ______, *On some mathematical principles in the linear theory of damped oscillations of continua I*, Integr. Equ. Oper. Theory **1** (1978), no. 3, 364–399.

[35] J. Labrousse, A. Sandovici, H. de Snoo, and H. Winkler, *The Kato decomposition of quasi-Fredholm relations*, Operators and Matrices **4** (2010), no. 1, 1–51.

[36] J. P. Labrousse, A. Sandovici, H. de Snoo, and H. Winkler, *Closed linear relations and their regular points*, Technical Report Preprint No. M 11/03, Technische Universität Ilmenau, Leiter des Instituts für Mathematik, 2011.

[37] H. Langer, *Invariant subspaces for a class of operators in spaces with indefinite metric*, J. Func. Anal. **19** (1975), no. 3, 232–241.

[38] S. Ôta, *Closed linear operators with domain containing their range*, Proceedings of the Edinburgh Math. Soc. **27** (1984), no. 2, 229–233.

[39] L. Pontryagin, *Hermitian operators in spaces with indefinite metric*, Izv. Akad. Nauk (in Russian) **8** (1944), no. 6, 243–280.

[40] D. Popovici and Z. Sebestyén, *Factorizations of linear relations*, Advances in Mathematics **233** (2013), 40–55.

[41] J. Rios-Cangas and L. O. Silva, *Dissipative extension theory for linear relations*, Expo. Math. **38** (2020), no. 1, 60–90.

[42] A. Sandovici, *On the adjoint of linear relations in Hilbert spaces*, Mediterr. J. Math. **17** (2020), no. 2, 1–23.

[43] A. Sandovici, H. de Snoo, and H. Winkler, *Ascent, descent, nullity, defect, and related notions for linear relations in linear spaces*, Lin. Alg. Appl. **423** (2007), no. 2–3, 456–497.

[44] K. Schmüdgen, *Unbounded Self-adjoint Operators on Hilbert Space*, Springer Dordrecht Heidelberg New York London, 2012.

[45] A. A. Shkalikov, *On invariant subspaces of dissipative operators in a space with indefinite metric*, arXiv:math (2004).

[46] P. Sorjonen, *On linear relations in an indefinite inner product space*, Annales Academire Scientiarum Fennicae **4** (1979), 169–192.

[47] V. A. Strauss and C. Trunk, *Some Sobolev spaces as Pontryagin spaces*, Vestn. Yuzhno-Ural. Gos. Un-ta. Ser. Matem. Mekh. Fiz. **6** (2012), 14–23.

[48] Z. Tarcsay and Z. Sebestyén, *Range-kernel characterizations of operators which are adjoint of each other*, Advances in Operator Theory **5** (2020), no. 3, 1026–1038.

VILNIUS UNIVERSITY, INSTITUTE OF THEORETICAL PHYSICS AND ASTRONOMY, SAULĖTEKIO AVE. 3, 10257 VILNIUS, LITHUANIA

*Email address*: `rytis.jursenas@tfai.vu.lt`